\documentclass[a4paper,twoside,11pt,english,intlimits]{article}

\usepackage{amsmath,amsthm,amsfonts,amssymb}
\usepackage{babel}
\usepackage[T1]{fontenc}
\usepackage{times,euler,euscript,stmaryrd,upgreek}
\usepackage{fancyhdr,titlesec,url,enumerate}
\usepackage[all]{xy}
\usepackage{hyperref}
\usepackage{catex}

\CompileMatrices

\DeclareFontFamily{OT1}{pzc}{}
\DeclareFontShape{OT1}{pzc}{m}{it}{<-> s * [1.2] pzcmi7t}{}
\DeclareMathAlphabet{\mathscr}{OT1}{pzc}{m}{it}
\titleformat{\section}[block]{\scshape\filcenter\Large}{\thesection.}{.5em}{}
\titleformat{\subsection}[block]{\bfseries\filcenter\large}{\thesubsection.}{.5em}{\medskip}
\titleformat{\subsubsection}[runin]{\bfseries}{\thesubsubsection.}{.5em}{}[.]

\newtheoremstyle{ntheorem}%
	{\topsep}{\topsep}{\itshape}{0pt}{\bfseries}{.}{.5em}%
	{\thmnumber{#2.\hspace{.5em}}\thmname{#1}\thmnote{ (#3)}}
	
\newtheoremstyle{ndefinition}%
	{\topsep}{\topsep}{\normalfont}{0pt}{\bfseries}{.}{.5em}%
	{\thmnumber{#2.\hspace{.5em}}\thmname{#1}\thmnote{ (#3)}}
	
\newtheoremstyle{nremark}%
	{\topsep}{\topsep}{\normalfont}{0pt}{\itshape}{.}{.5em}%
	{\thmnumber{}\thmname{#1}\thmnote{ (#3)}}

\theoremstyle{ntheorem}
  	\newtheorem{theorem}[subsubsection]{Theorem}
  	\newtheorem{proposition}[subsubsection]{Proposition}
	\newtheorem{lemma}[subsubsection]{Lemma}
  	\newtheorem{corollary}[subsubsection]{Corollary}

\theoremstyle{ndefinition}

	\newtheorem{example}[subsubsection]{Example}

\theoremstyle{nremark}

\fancypagestyle{plain}{
\fancyhf{}\fancyfoot[LC,RC]{}
\fancyfoot[LE,RO]{$\thepage$}

}

\UseTips
\SelectTips{eu}{11}

\newdir{ >}{{}*!/-10pt/@{>}}
\newdir{ -}{{}*!/-10pt/@{}}
\newdir{> }{{}*!/+10pt/@{>}}

\makeatletter

\xyletcsnamecsname@{dir4{}}{dir{}}
\xydefcsname@{dir4{-}}{\line@ \quadruple@\xydashh@}
\xydefcsname@{dir4{.}}{\point@ \quadruple@\xydashh@}
\xydefcsname@{dir4{~}}{\squiggle@ \quadruple@\xybsqlh@}
\xydefcsname@{dir4{>}}{\Tttip@}
\xydefcsname@{dir4{<}}{\reverseDirection@\Tttip@}

\xydef@\quadruple@#1{%
	\edef\Drop@@{%
		\dimen@=#1\relax
		\dimen@=.5\dimen@
		\A@=-\sinDirection\dimen@
		\B@=\cosDirection\dimen@
		\setboxz@h{%
			\setbox2=\hbox{\kern3\A@\raise3\B@\copy\z@}%
			\dp2=\z@ \ht2=\z@ \wd2=\z@ \box2
			\setbox2=\hbox{\kern\A@\raise\B@\copy\z@}%
			\dp2=\z@ \ht2=\z@ \wd2=\z@ \box2
			\setbox2=\hbox{\kern-\A@\raise-\B@\copy\z@}%
			\dp2=\z@ \ht2=\z@ \wd2=\z@ \box2
			\setbox2=\hbox{\kern-3\A@\raise-3\B@ \noexpand\boxz@}%
			\dp2=\z@ \ht2=\z@ \wd2=\z@ \box2
		}%
		\ht\z@=\z@ \dp\z@=\z@ \wd\z@=\z@ \noexpand\styledboxz@
	}%
}

\xydef@\Tttip@{\kern2pt \vrule height2pt depth2pt width\z@
	\Tttip@@ \kern2pt \egroup
	\U@c=0pt \D@c=0pt \L@c=0pt \R@c=0pt \Edge@c={\circleEdge}%
	\def\Leftness@{.5}\def\Upness@{.5}%
	\def\Drop@@{\styledboxz@}\def\Connect@@{\straight@{\dottedSpread@\jot}}}
	
\xydef@\Tttip@@{%
	\dimen@=.25\dimen@
 	\B@=\cosDirection\dimen@
	\setboxz@h\bgroup\reverseDirection@\line@ \wdz@=\z@ \ht\z@=\z@ \dp\z@=\z@
	{\vDirection@(1,-1)\xydashl@ \xyatipfont\char\DirectionChar}%
	{\vDirection@(1,+1)\xydashl@ \xybtipfont\char\DirectionChar}%
}

\xydef@\ar@form{
	\ifx \space@\next \expandafter\DN@\space{\xyFN@\ar@form}%
	\else\ifx ^\next \DN@ ^{\xyFN@\ar@style}\edef\arvariant@@{\string^}%
	\else\ifx _\next \DN@ _{\xyFN@\ar@style}\edef\arvariant@@{\string_}%
	\else\ifx 0\next \DN@ 0{\xyFN@\ar@style}\def\arvariant@@{0}%
	\else\ifx 1\next \DN@ 1{\xyFN@\ar@style}\def\arvariant@@{1}%
	\else\ifx 2\next \DN@ 2{\xyFN@\ar@style}\def\arvariant@@{2}%
	\else\ifx 3\next \DN@ 3{\xyFN@\ar@style}\def\arvariant@@{3}%
	\else\ifx 4\next \DN@ 4{\xyFN@\ar@style}\def\arvariant@@{4}%
	\else\ifx \bgroup\next \let\next@=\ar@style
	\else\ifx [\next \DN@[##1]{\ar@modifiers{[##1]}}
	\else\ifx *\next \DN@ *{\ar@modifiers}%
	\else\addLT@\ifx\next \let\next@=\ar@slide
	\else\ifx /\next \let\next@=\ar@curveslash
	\else\ifx (\next \let\next@=\ar@curveinout 
	\else\addRQ@\ifx\next \addRQ@\DN@{\ar@curve@}%
	\else\addLQ@\ifx\next \addLQ@\DN@{\xyFN@\ar@curve}%
	\else\addDASH@\ifx\next \addDASH@\DN@{\defarstem@-\xyFN@\ar@}%
	\else\addEQ@\ifx\next \addEQ@\DN@{\def\arvariant@@{2}\defarstem@-\xyFN@\ar@}%
	\else\addDOT@\ifx\next \addDOT@\DN@{\defarstem@.\xyFN@\ar@}%
	\else\ifx :\next \DN@:{\def\arvariant@@{2}\defarstem@.\xyFN@\ar@}%
	\else\ifx ~\next \DN@~{\defarstem@~\xyFN@\ar@}%
	\else\ifx !\next \DN@!{\dasharstem@\xyFN@\ar@}%
	\else\ifx ?\next \DN@?{\ar@upsidedown\xyFN@\ar@}%
	\else \let\next@=\ar@error
	\fi\fi\fi\fi\fi\fi\fi\fi\fi\fi\fi\fi\fi\fi\fi\fi\fi\fi\fi\fi\fi\fi\fi \next@}

\makeatother

\newcommand{\fl}{\to}

\newcommand{\dfl}{\Rightarrow}
\newcommand{\tfl}{\Rrightarrow}
\newcommand{\qfl}{\!\raisebox{0.3em}{\xymatrix@C=10pt{ \ar@4[r] &}}\!}

\newcommand{\ens}[1]{\left\{#1\right\}}
\newcommand{\iar}[1]{\lfloor#1\rfloor}

\newcommand{\cl}[1]{\overline{#1}}
\renewcommand{\tilde}[1]{\widetilde{#1}}
\newcommand{\rep}[1]{\widehat{#1}}

\newcommand{\tens}{\otimes}

\newcommand{\ie}{\emph{i.e.}}
\renewcommand{\phi}{\varphi}
\renewcommand{\epsilon}{\varepsilon}
\newcommand{\dr}{\partial}
\newcommand{\Nb}{\mathbb{N}}
\newcommand{\Ar}{\EuScript{A}}
\newcommand{\Br}{\EuScript{B}}
\newcommand{\Cr}{\EuScript{C}}
\newcommand{\Dr}{\EuScript{D}}
\renewcommand{\Pr}{\EuScript{P}}
\newcommand{\C}{\mathbf{C}}
\newcommand{\T}{\mathbf{T}}
\renewcommand{\P}{\mathbf{P}}
\newcommand{\Q}{\mathbf{Q}}
\newcommand{\poly}[1]{\mathrm{#1}}

\newcommand{\Cat}{\mathbf{Cat}}
\newcommand{\Perm}{\mathbf{Perm}}
\newcommand{\As}{\mathbf{As}}
\newcommand{\Mon}{\mathbf{Mon}}
\newcommand{\Sym}{\mathbf{Sym}}
\newcommand{\Brd}{\mathbf{Br}}

\newcommand{\Alg}{\mathscr{Alg}}
\newcommand{\ab}[1]{#1_{\text{ab}}}
\newcommand{\tck}[1]{{#1}^{\top}}
\newcommand{\abtck}[1]{\ab{\tck{#1}}}
\DeclareMathOperator{\Aut}{Aut}
\DeclareMathOperator{\id}{id}

\newcommand{\pdf}[1]{\texorpdfstring{$#1$}{1}}

\deftwocell[black]{mu : 2 -> 1}
\deftwocell[black,circle]{eta : 0 -> 1}
\deftwocell[crossing]{tau : 2 -> 2}
\deftwocell[gray]{alpha : 3 -> 1}
\deftwocell[gray,lefthalfcircle]{lambda : 1 -> 1}
\deftwocell[gray,righthalfcircle]{rho : 1 -> 1}
\deftwocell[gray,rectangle]{beta : 2 -> 1}
\deftwocell[gray,rectangle]{gamma : 3 -> 1}
\deftwocell[white]{aleph : 4 -> 1}
\deftwocell[white]{beth : 2 -> 1}
\deftwocell[white,rectangle]{gimmel : 2 -> 1}
\deftwocell[white,rectangle]{daleth : 3 -> 1}

\deftwocell[crossing1]{tau1 : 3 -> 3}
\deftwocell[crossing2]{tau2 : 3 -> 3}
\deftwocell[dots]{dots : 2 -> 2}
\deftwocell[gray,rectangle]{f : 2 -> 2}
\deftwocell[lightgray,rectangle]{g : 2 -> 2}
\deftwocell[gray,halfcircle]{phi : 2 -> 1}
\deftwocell[lightgray,halfcircle]{psi : 2 -> 1}

\begin{document}

\thispagestyle{empty}
\begin{center}
\textbf{\LARGE Coherence in monoidal track categories}

\vspace{4mm}
\begin{tabular}{c c c}
\textbf{\large Yves Guiraud} &$\qquad\qquad$& \textbf{\large Philippe Malbos} \\
INRIA && Université Lyon 1 \\
Institut Camille Jordan && Institut Camille Jordan \\
guiraud@math.univ-lyon1.fr && malbos@math.univ-lyon1.fr
\end{tabular}
\end{center} 

\vspace{2mm}
\begin{em}
\hrule height 1.5pt

\medskip
\noindent \textbf{Abstract --} We introduce homotopical methods based on rewriting on higher-dimensional categories to prove coherence results in categories with an algebraic structure. We express the coherence problem for (symmetric) monoidal categories as an asphericity problem for a track category and we use rewriting methods on polygraphs to solve it. The setting is extended to more general coherence problems, seen as 3-dimensional word problems in a track category, including the case of braided monoidal categories.


\noindent
\textbf{M.S.C. 2000 --} 18C10, 18D10, 68Q42.

\medskip
\hrule height 1.5pt
\end{em}

\section*{Introduction}

A monoidal category is a category equipped with a product, associative up to a natural isomorphism, and having a distinguished object, which is a unit for the product up to natural isomorphisms. Associativity and unity satisfy, in turn, a coherence condition: all the diagrams built from the corresponding natural isomorphisms are commutative. A cornerstone result for monoidal categories was to reduce the infinite requirement ``every diagram commutes'' to a finite requirement ``if a specified finite set of diagrams commute then every diagram commutes'',~\cite{MacLane63, Stasheff63}. We call \emph{coherence basis} such a finite set of diagrams. 

A symmetric monoidal category is a monoidal category whose product is commutative up to a natural isomorphism, called symmetry. In a symmetric monoidal category the coherence problem has the same formulation as in monoidal categories, with additional coherence diagrams for the symmetry,~\cite{MacLane63}.

In a symmetric monoidal category the symmetry is its own inverse. Braided monoidal categories are monoidal categories commutative up to an isomorphism which is not its own inverse. The coherence problem in braided categories has another formulation: a diagram is commutative if and only if its two sides correspond to the same braid,~\cite{JoyalStreet93}.

In this paper, we formulate the coherence problem for monoidal track $2$-categories in the homotopical terms of higher-dimensional categories, as introduced by the authors in~\cite{GuiraudMalbos09}. This formulation gives a way to reduce the coherence problem to a $3$-dimensional word problem in track categories. The construction of convergent (\ie, terminating and confluent) presentations of monoidal track $2$-categories allows us to reduce the problem ``every diagram commutes'' to ``if the diagrams induced by critical branchings commute then every diagram commutes'': the confluence diagrams of critical branchings form a coherence basis. Let us illustrate this methodology on a simple example.

\subsection*{Coherence for categories with an associative product}

Let us consider a category $\Cr$ equipped with a functor $\tens:\Cr\times\Cr\fl\Cr$ which is associative up to a natural isomorphism, \ie, there is a natural isomorphism 
\[
\alpha_{x,y,z} \: : \: (x\tens y)\tens z \:\longrightarrow\: x\tens(y\tens z) \:,
\]
such that the following diagram commutes in $\Cr$:
\begin{equation}
\label{diagramPenta}
\scalebox{0.8}{
\xymatrix@C=-1em{
& { (x\tens (y\tens z))\tens t}
	\ar[rr] ^-{\alpha}
&& { x\tens ((y\tens z) \tens t)}
	\ar[dr] ^-{\alpha}
\\
{ ((x\tens y)\tens z) \tens t}
	\ar[ur] ^-{\alpha}
	\ar[drr] _-{\alpha}
&& { \copyright}
&& { x\tens(y\tens (z\tens t))}
\\
&& { (x\tens y)\tens (z\tens t)}
	\ar[urr] _-{\alpha}
}
}
\end{equation}
Presentation of such categories by generators and relations can be achieved using the notion of   polygraph. This notion of presentation of higher-dimensional categories was introduced by Burroni, \cite{Burroni93}, and by Street under the terminology of computads,~\cite{Street76,Street87}. In this paper, we use Burroni's terminology, as usual in rewriting theory. An \emph{$n$-polygraph} is a family $(\Sigma_0,\ldots,\Sigma_n)$, where $\Sigma_0$ is a set and, for every $0\leq k < n$, $\Sigma_{k+1}$ is a family of parallel $k$-cells of the free $k$-category   $\Sigma_k^\ast$ over $\Sigma_k$. We call such a family a \emph{cellular extension of $\Sigma_k^\ast$}.
  
Categories with an associative product can be presented using the notion of polygraph as follows.
Let us consider the $3$-polygraph $\poly{As}_3$ with one $0$-cell, one $1$-cell $\;\twocell{1}\;$, one $2$-cell $\twocell{mu}$ and one $3$-cell: 
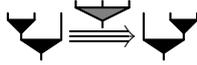
\[
\xymatrix{
{\twocell{ (mu *0 1) *1 mu }}
	\ar@3[r] ^-{\twocell{alpha}}
& {\twocell{ (1 *0 mu) *1 mu } }
}
\]
Let $\tck{\poly{As}}_3$ be the free \emph{track $3$-category} generated by $\poly{As}_3$, \ie, the free $3$-category over $\poly{As}_3$ whose $3$-cells are invertible. The relation~\eqref{diagramPenta} satisfied by the associativity isomorphism can be presented by a cellular extension $\poly{As}_4$ of
$\tck{\poly{As}}_3$ with one $4$-cell:
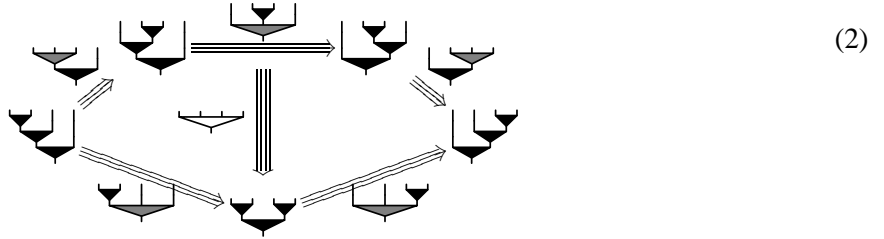
\begin{equation}
\label{4cellAleph}
\scalebox{0.8}{
\xymatrix@R=10pt@C=15pt{
&{ \twocell{ (1 *0 mu *0 1) *1 (mu *0 1) *1 mu } }
	\ar@3 [rr] 
		^-{\twocell{(1 *0 mu *0 1) *1 alpha}} 
		_-{}="1" 
&& { \twocell{ (1 *0 mu *0 1) *1 (1 *0 mu) *1 mu } }
	\ar@3 [dr] ^-{\twocell{(1 *0 alpha) *1 mu} }
\\
{ \twocell{ (mu *0 2) *1 (mu *0 1) *1 mu } }
	\ar@3 [ur] ^-{\twocell{(alpha *0 1) *1 mu}}
	\ar@3 [drr] _-{\twocell{(mu *0 2)*1 alpha}}
&&&& { \twocell{ (2 *0 mu) *1 (1 *0 mu) *1 mu } }
\\
&& { \twocell{ (mu *0 mu) *1 mu } }
	\ar@3 [urr] _-{\twocell{(2 *0 mu)*1 alpha}}
	\ar@4 "1"!<0pt,-10pt>;[]!<0pt,20pt> _-*+{\twocell{aleph}\;}
}
}
\end{equation}
Let $\As\Cat$ be the track $3$-category obtained as the quotient of $\tck{\poly{As}}_3$ by the cellular extension~$\poly{As}_4$. The category of (small) categories with a product, associative up to a natural isomorphism, is isomorphic to the category $\Alg(\As\Cat)$ of algebras over the $3$-category $\As\Cat$. Such an algebra is a $3$-functor from $\As\Cat$ to the monoidal $2$-category $\Cat$ of small categories, functors and natural transformations, seen as a $3$-category with only one $0$-cell. The correspondence associates, to a category $(\Cr,\tens,\alpha)$, the algebra $\Ar : \As\Cat \fl \Cat$ defined by:
\[
\Ar(\:\twocell{1}\:) \:=\: \Cr,
\qquad\quad
\Ar(\twocell{mu}) \:=\: \tens,
\quad\qquad
\Ar(\twocell{alpha}) \:=\: \alpha.
\]
A diagram in a $\As\Cat$-algebra $\Ar$ is the image $\Ar(\gamma)$ of a pair $\gamma=(A,B)$ of parallel $3$-cells in $\As\Cat$. This diagram commutes if $\Ar(A)=\Ar(B)$ holds in $\Cat$. The coherence problem for $\As\Cat$-algebras can be formulated as ``does every diagram commute in every $\As\Cat$-algebra''. In this way, the coherence problem is reduced to showing that~$\poly{As}_4$ forms a \emph{homotopy basis} of the track $3$-category $\tck{\poly{As}}_3$, \ie, if $(A,B)$ is a pair of parallel $3$-cells of $\As\Cat$, then $A=B$. A $3$-category that satisfies this last property is called \emph{aspherical}.

Proving asphericity is a special case of a word problem in a track $3$-category. In~\cite{GuiraudMalbos09}, the authors prove that for a convergent, \ie, terminating and confluent, $n$-polygraph $\Sigma$, the \emph{critical branchings} generate a homotopy basis of the free track $n$-category $\tck{\Sigma}$. In our example, the $3$-polygraph $\poly{As}_3$ is convergent and has a unique critical branching, formed by the two different applications of the $3$-cell~\twocell{alpha} on the same $2$-cell:
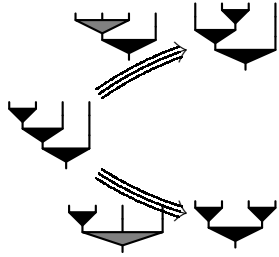
\[
\xymatrix@R=0.5em@C=3em{
&{ \twocell{ (1 *0 mu *0 1) *1 (mu *0 1) *1 mu } }
\\
{ \twocell{ (mu *0 2) *1 (mu *0 1) *1 mu } }
	\ar@/^/@3 [ur] ^-{\twocell{(alpha *0 1) *1 mu}}
	\ar@/_/@3 [dr] _-{\twocell{(mu *0 2)*1 alpha}}
\\
& { \twocell{ (mu *0 mu) *1 mu } }
}
\]
The $4$-cell \twocell{aleph} of~\eqref{4cellAleph} forms a \emph{confluence diagram} for this critical branching. As a consequence, the cellular extension $\poly{As}_4$ is a homotopy basis of the track $3$-category $\tck{\poly{As}}_3$: this proves the coherence result for $\As\Cat$-algebras. 

\subsection*{Organisation of the paper}

In Section~\ref{Section:Preliminaries}, we recall notions on higher-dimensional track categories, presentations by polygraphs and polygraphic rewriting, including critical branchings. We introduce the notion of higher-dimensional pro(p)s in Section~\ref{Subsection:HigherDimensionalPRO(P)S}. 

For $n\geq 1$, a \emph{(track) $n$-pro} is a (track) $n$-category with one $0$-cell, such that its underlying $1$-category is the monoid of natural numbers with addition. Equivalently, for $n\geq 2$, a (track) $n$-pro is a strict monoidal category (seen as a $2$-category with one $0$-cell), enriched in (track) $(n-2)$-categories and whose underlying monoid of objects is the monoid of natural numbers with the addition. A \emph{(track) $n$-prop} is a (track) $n$-pro, whose underlying monoidal category is symmetric. In particular, $2$-pro(p)s coincide with Mac Lane's PRO(P)s, an acronym for ``product (and permutation) categories", introduced in~\cite{MacLane65}. 

For coherence problems, we consider special cases of track $3$-pro(p)s: the track $3$-pros $\As\Cat$ of categories with an associative product and $\Mon\Cat$ of monoidal categories and the track $3$-props $\Sym\Cat$ of symmetric monoidal categories and $\Brd\Cat$ of braided monoidal categories.

An \emph{algebra over a $3$-pro(p)} $\P$ is a strict (symmetric) monoidal $2$-functor from $\P$ to $\Cat$. Here $\Cat$ is considered as a $3$-category with one $0$-cell, categories as $1$-cells, functors as $2$-cells and natural transformations as $3$-cells, see Paragraph~\ref{algebrasOverPro(p)}. In Proposition~\ref{Propaspherical=>coherence}, we relate the coherence problem for algebras over a $3$-pro(p) $\P$ to the asphericity of $\P$: if the $3$-pro(p) $\P$ is aspherical, then every $\P$-diagram commutes in every $\P$-algebra.

Thus, reducing the coherence problem ``every diagram commutes'' to ``if some diagrams commute then every diagram commutes'' consists in constructing an algebraic presentation of the $3$-pro(p) proving that it is aspherical. We show that a convergent presentation gives a procedure to solve the coherence problem.

\subsubsection*{The monoidal coherence problem} 

In Section~\ref{Section:CoherenceMonoidalCategories}, we consider the case of $3$-pros. A convergent presentation for a $3$-pro $\P$ is a pair $(\Sigma_3,\Sigma_4)$, where $\Sigma_3$ is a convergent $3$-polygraph together with a cellular extension $\Sigma_4$ of generating confluences of $\Sigma_3$. We have:

\medskip
\noindent 
\textbf{Theorem~\ref{TheoremProAspherical}.}
\emph{If a track $3$-pro $\P$ admits a convergent presentation, then every $\P$-diagram commutes in every $\P$-algebra.}

\medskip
\noindent 
In Section~\ref{Subsection:ApplicationCoherenceMon}, we consider the coherence problem for monoidal categories. We prove that the $3$-pro $\Mon\Cat$ of monoidal categories is aspherical, see Paragraph~\ref{corollaryMonAspherical}, hence the coherence theorem for monoidal categories, proved in~\cite{MacLane63}.

\subsubsection*{The symmetric monoidal coherence problem} 

In Section~\ref{Section:CoherenceSym}, for the coherence problem for symmetric monoidal categories, we consider the asphericity problem of \emph{algebraic} track $3$-props, \ie, track $3$-props whose generating $2$-cells and $3$-cells have coarity~$1$, see Paragraph~\ref{algebraicPresentationProp}. In that case, we have a convergent presentation of the symmetry, see~\cite{Burroni93,Guiraud04}. This gives the following sufficient condition for proving that an algebraic track $3$-prop is aspherical, where $\pi(\Gamma_{\Sigma_3})$ is a cellular extension generated by the critical branchings that do not dependent on the symmetry only: 

\medskip
\noindent 
\textbf{Theorem~\ref{TheoremPropAspherical}.}
\emph{If a track $3$-prop $\P$ admits an algebraic convergent presentation $(\Sigma_3,\Sigma_4)$ such that $\Sigma_4$ is Tietze-equivalent to $\pi(\Gamma_{\Sigma_3})$, then $\P$ is aspherical.}

\medskip
\noindent 
In the case of the $3$-prop $\Sym\Cat$ of symmetric monoidal categories, this result gives the corresponding coherence theorem, see Corollary~\ref{CoherenceTheoremSymmetricMonoidalCategories}.

\subsubsection*{The braided monoidal case and the generalised coherence problem} 

For braided monoidal categories, we consider a generalised version of the coherence problem: ``given a $2$-prop $\P$, decide, for any $3$-sphere $\gamma$ of $\P$, whether or not the diagram $\Ar(\gamma)$ commutes in every $\P$-algebra $\Ar$''. To solve it, we proceed in two steps. First, we prove that coherence is preserved by aspherical quotients, so that we can reduce a $3$-prop to its non-aspherical part: 

\medskip
\noindent 
\textbf{Theorem~\ref{thmEquivProps}.}
\emph{Let $\P$ and $\Q$ be $3$-props with $\Q$ aspherical and $\Q\subseteq\P$. Then, for every $3$-sphere $(A,B)$ of $\P$, we have $A=B$ if and only if $\pi(A)=\pi(B)$.} 

\medskip
\noindent 
Then, given an algebraic $3$-prop $\P$, we define the \emph{initial $\P$-algebra} $\Pr$, see Section~\ref{initialAlgebra}, and we prove: 

\medskip
\noindent 
\textbf{Theorem~\ref{theoremGeneCoherenceProblemAlg2Prop}.}
\emph{Let $\P$ be an algebraic $3$-prop and let $(A,B)$ be a $3$-sphere of $\P$. Then we have $A=B$ if and only if $\Pr(A)=\Pr(B)$.}

\medskip
\noindent 
In the case of the $3$-prop of braided monoidal categories, the initial algebra $\Br$ associates, to every $3$-cell $A$, a braid $\Br(A)$. Hence, the introduced methodology recovers the coherence result of Joyal and Street,~\cite{JoyalStreet93}: a diagram commutes if and only if its two sides are associated to the same braid.

\section{Preliminaries}
\label{Section:Preliminaries}

In this section, we recall from~\cite{GuiraudMalbos09} notions and results on higher-dimensional (track) categories, homotopy bases and presentations by polygraphs.

\subsection{Higher-dimensional categories and homotopy bases}

Let $n$ be a natural number and let $\C$ be an $n$-category (we always consider strict, globular $n$-categories). We denote by $\C_k$ the set (and the $k$-category) of $k$-cells of~$\C$. If $f$ is in~$\C_k$, then $s_i(f)$ and $t_i(f)$ respectively denote the $i$-source and $i$-target of $f$; we drop the suffix~$i$ when $i=k-1$. The source and target maps satisfy the \emph{globular relations}: 
\[
s_i\circ s_{i+1} \:=\: s_i\circ t_{i+1} 
\qquad\text{and}\qquad
t_i\circ s_{i+1} \:=\: t_i\circ t_{i+1}.
\]
We respectively denote by $f:u\fl v$, $\;f:u\dfl v$, $\;f:u\tfl v\;$ or $\;f:u\qfl v\;$ a $1$-cell, $2$-cell, $3$-cell or $4$-cell $f$ with source~$u$ and target~$v$. 

If $(f,g)$ is a pair of $i$-composable $k$-cells, that is when $t_i(f)=s_i(g)$, we denote by $f\star_i g$ their $i$-composite. The compositions satisfy the \emph{exchange relations} given, for every $i\neq j$ and every possible cells $f$, $g$, $h$ and $k$, by: 
\[
(f \star_i g) \star_j (h \star_i k) \:=\: (f \star_j h) \star_i (g \star_j k).
\]
If $f$ is a $k$-cell, we denote by $\id_f$ its identity $(k+1)$-cell. When $\id_f$ is composed with cells of dimension $k+1$ or higher, we simply denote it by~$f$. A cell is \emph{degenerate} when it is an identity cell.

\subsubsection{Track \pdf{n}-categories}

In an $n$-category~$\C$, we say that a $k$-cell $f$ with source $u$ and target $v$ is \emph{invertible} when it admits an inverse for the higher-dimensional composition $\star_{k-1}$ defined on it, \ie, when there exists a (necessarily unique) $k$-cell in $\C$, with source $v$ and target $u$ in $\C$, denoted by $f^-$ and called the \emph{inverse of $f$}, that satisfies 
\[
f\star_{k-1} f^- \:=\: \id_u 
\qquad\text{and}\qquad
f^-\star_{k-1} f \:=\: \id_v.
\]
A \emph{track $n$-category} is an $n$-category whose $n$-cells are invertible. One can also define track $n$-categories by induction on $n$, with track $1$-categories being groupoids and track $(n+1)$-categories being categories enriched in track $n$-categories.

\subsubsection{Cellular extensions}

Let $\C$ be an $n$-category. A \emph{$k$-sphere of $\C$} is a pair $\gamma=(f,g)$ of parallel $k$-cells of $\C$, \ie, with $s(f)=s(g)$ and $t(f)=t(g)$. We call $f$ the \emph{source} of $\gamma$ and $g$ its \emph{target}. When $f=g$, the $k$-sphere $\gamma$ is \emph{degenerate}. An $n$-category $\C$ is \emph{aspherical} when every $n$-sphere of $\C$ is degenerate.

A \emph{cellular extension of $\C$} is a family $\Gamma$ of $n$-spheres of $\C$. By considering all the formal compositions of elements of $\Gamma$, seen as $(n+1)$-cells with source and target in $\C$, one builds the \emph{free $(n+1)$-category generated by $\Gamma$ over $\C$}, denoted by~$\C[\Gamma]$. 

The \emph{quotient of $\C$ by $\Gamma$}, denoted by $\C/\Gamma$, is the $n$-category one gets from $\C$ by identification of the $n$-cells $s(\gamma)$ and $t(\gamma)$ for every element $\gamma$ of $\Gamma$. Two cellular extensions $\Gamma_1$ and $\Gamma_2$ of $\C$ are \emph{Tietze-equivalent} if the $n$-categories $\C/\Gamma_1$ and $\C/\Gamma_2$ are isomorphic. 

The \emph{free track $(n+1)$-category generated by $\Gamma$ over $\C$} is defined by
\[
\C(\Gamma) \:=\: \C[\Gamma,\Gamma^-]/\mathrm{Inv}(\Gamma),
\]
where $\Gamma^-$ and $\mathrm{Inv}(\Gamma)$ are the following cellular extensions of $\C$ and $\C[\Gamma,\Gamma^-]$, respectively:
\[
\Gamma^- \:=\: \ens{ \:\gamma^-: t(\gamma) \fl s(\gamma)
  \;\;|\;\; \gamma\in\Gamma \:} 
\]
and
\[
\mathrm{Inv}(\Gamma) \:=\: 
	\ens{ \: \gamma\star_n\gamma^- \fl \id_{s\gamma} \;,\: 
	\gamma^-\star_n\gamma \fl \id_{t\gamma} \;\;|\;\; 
	\gamma\in \Gamma }.
\]

\subsubsection{Homotopy bases}

A cellular extension $\Gamma$ of an $n$-category $\C$ is a \emph{homotopy basis} when the quotient $n$-category $\C/\Gamma$ is aspherical, \ie, when, for every $n$-sphere $\gamma$ of $\C$, there exists an $(n+1)$-cell from $s(\gamma)$ to $t(\gamma)$ in the track $(n+1)$-category $\C(\Gamma)$.

\subsection{Presentations by polygraphs}
\label{Polygraphs}

We define, by induction on $n$, the notions of $n$-polygraph, of presented $(n-1)$-category and of freely generated (track) $n$-category. For a deeper treatment, we refer the reader to~\cite{Burroni93,Metayer03,GuiraudMalbos09}.

A \emph{$1$-polygraph} is a graph $\Sigma=(\Sigma_0,\Sigma_1)$. We denote by $\Sigma^*$ the free $1$-category and by $\tck{\Sigma}$ the free track $1$-category (\ie, groupoid) it generates. An $(n+1)$-polygraph is a pair $\Sigma=(\Sigma_n,\Sigma_{n+1})$ made of an $n$-polygraph $\Sigma_n$ and a cellular extension~$\Sigma_{n+1}$ of the free $n$-category $\Sigma_n^*$ generated by the $n$-polygraph~$\Sigma_n$. The \emph{$n$-category presented by $\Sigma$}, the \emph{free $(n+1)$-category generated by $\Sigma$} and the \emph{free track $(n+1)$-category generated by $\Sigma$} are respectively denoted by $\cl{\Sigma}$, $\Sigma^*$ and $\tck{\Sigma}$ and defined as follows: 
\[
\cl{\Sigma} \:=\: \Sigma_n^*/\Sigma_{n+1},
\qquad\qquad
\Sigma^* \:=\: \Sigma_n^*[\Sigma_{n+1}],
\qquad\qquad
\tck{\Sigma} \:=\: \Sigma_n^*(\Sigma_{n+1}).
\]
An $n$-polygraph yields a diagram of cellular extensions, as given in~\cite{Burroni93}:
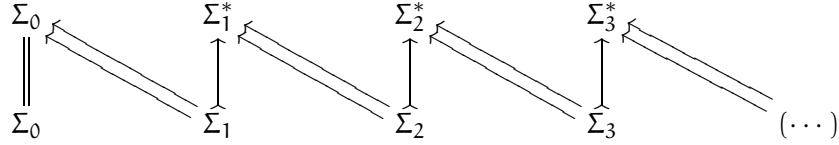
\[
\xymatrix{
\Sigma_0
     \ar@{=}[d] 
&& \Sigma_1^* 
&& \Sigma_2^*
&& \Sigma_3^*
\\
\Sigma_0
&& \Sigma_1
     \ar@<0.5ex>[llu]
     \ar@<-0.5ex>[llu]
     \ar@{>->}[u]
&& \Sigma_2
     \ar@<0.5ex>[llu]
     \ar@<-0.5ex>[llu]
     \ar@{>->}[u]
&& \Sigma_3
     \ar@<0.5ex>[llu]
     \ar@<-0.5ex>[llu]
     \ar@{>->}[u]
&& (\cdots)
     \ar@<0.5ex>[llu]
     \ar@<-0.5ex>[llu]
}
\]
If $\C$ is an $n$-category, a \emph{presentation of $\C$} is an $(n+1)$-polygraph $\Sigma$ such that $\cl{\Sigma}$ is isomorphic to $\C$.

\subsubsection{Polygraphic rewriting}

Let $\Sigma$ be an $n$-polygraph. We say that an $(n-1)$-cell $u$ of $\Sigma$ \emph{reduces} to some $(n-1)$-cell $v$ in $\Sigma$ when there exists a non-degenerate $n$-cell from $u$ to $v$ in $\Sigma^*$. A \emph{reduction sequence of $\Sigma$} is a countable family $(u_i)_{i\in I}$ of $(n-1)$-cells of $\Sigma$ such that each $u_i$ reduces to the follo\-wing~$u_{i+1}$. We say that $\Sigma$ \emph{terminates} when it has no infinite reduction sequence.

A \emph{branching} of $\Sigma$ is a non-ordered pair $(f,g)$ of $n$-cells of $\Sigma^*$ with the same source, called the \emph{source of $(f,g)$}. A branching $(f,g)$ is \emph{confluent} when there exists a pair $(f',g')$ of $n$-cells of $\Sigma^*$ with the same target and such that $(f,f')$ and $(g,g')$ are composable, as in the following diagram, called a \emph{confluence diagram for the branching $(f,g)$}:
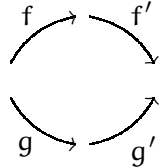
\[
\xymatrix@R=1em@W=0.5em@H=0.5em@M=0em{
& \strut
	\ar@/^/[dr] ^-{f'} 
\\
\strut 
	\ar@/^/[ur] ^-{f}
	\ar@/_/[dr] _-{g}
&& \strut
\\
& \strut
	\ar@/_/[ur] _-{g'}
}
\]
We say that the $n$-polygraph $\Sigma$ is \emph{confluent} when every branching of $\Sigma$ is confluent.
 
Finally, the $n$-polygraph $\Sigma$ is \emph{convergent} when it terminates and it is confluent. Following~\cite{Thue14}, finite and convergent rewriting systems, such as convergent polygraphs, give an algorithmic way, the \emph{normal form algorithm}, to solve the word problem for the algebraic structure they present: see~\cite{BookOtto93} for presentations of monoids by word (or string) rewriting systems,~\cite{BaaderNipkow98} for presentations of equational theories by term rewriting systems and~\cite{GuiraudMalbos09} for presentations of $n$-categories by polygraphs. Here, we are interested in convergent polygraphs because they give a way to compute homotopy bases.

\subsubsection{Critical branchings}

Here, we give the informal idea underlying the notion of critical branchings. We refer the reader to other works for a fuller treatment of the subject:~\cite{BookOtto93} for word rewriting systems;~\cite{BaaderNipkow98} for term rewriting systems;~\cite{GuiraudMalbos09}, where the authors give a general theory of branchings in $n$-polygraphs and a thorough study of critical branchings of $3$-polygraphs;~\cite{ktheory}, where the authors describe resolutions of small categories based on the critical branchings (and generalisations) of presentations by convergent $2$-polygraphs.

Branchings in an $n$-polygraph $\Sigma$ occur when an $(n-1)$-cell $u$ of $\Sigma^*$ contains the sources of two $n$-cells $\phi$ and $\psi$ of $\Sigma$. When those sources are disjoint in $u$, the branching is confluent, such as in the following simple case, with $0\leq i\leq n-1$:
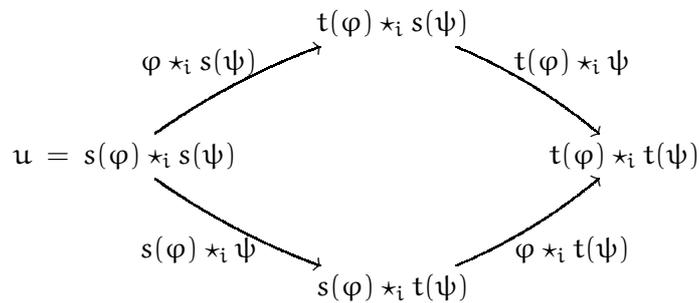
\[
\xymatrix@R=3em{
& t(\phi)\star_i s(\psi)
	\ar@/^/ [dr] ^-{t(\phi)\star_i \psi}
\\
u\:=\: s(\phi) \star_i s(\psi)
	\ar@/^/ [ur] ^-{\phi\star_i s(\psi)}
	\ar@/_/ [dr] _-{s(\phi)\star_i\psi}
&& t(\phi)\star_i t(\psi)
\\
& s(\phi) \star_i t(\psi)
	\ar@/_/ [ur] _-{\phi\star_i t(\psi)}
}
\]
Note that, in this example, both composites are equal to the $n$-cell $\phi\star_i \psi$, due to the exchange relation between $\star_i$ and $\star_n$ in $\Sigma^*$.

Otherwise, when the sources of $\phi$ and $\psi$ overlap in $u$, in such a way that $u$ is a minimal $(n-1)$-cell such that this overlapping occurs, we have a \emph{critical} branching. For example, in  a $2$-polygraph, we can have two different shapes of critical branchings:
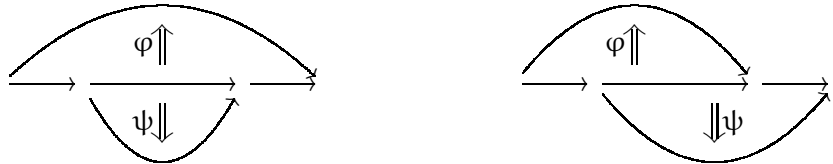
\[
\xymatrix{
\strut 
	\ar[r] _-{}
	\ar@/^6ex/ [rrrr] _-{}="t1"
& \strut
	\ar[rr] |-{} ^-{}="s1" _-{}="s2"
	\ar@/_6ex/ [rr] ^-{}="t2"
&& \strut
	\ar[r] _-{}
& \strut
\ar@2 "s1"!<0pt,7.5pt>;"t1"!<0pt,-7.5pt> ^-{\phi}
\ar@2 "s2"!<0pt,-7.5pt>;"t2"!<0pt,7.5pt> _-{\psi}
}
\qquad\qquad\qquad
\xymatrix{
\strut 
	\ar[r] _-{}
	\ar@/^6ex/ [rrr] _-{}="t1"
& \strut
	\ar[rr] |-{} _(0.75){}="s2" ^(0.25){}="s1"
	\ar@/_6ex/ [rrr] ^-{}="t2"
&& \strut
	\ar[r] ^-{}
& \strut
\ar@2 "s1"!<0pt,7.5pt>;"t1"!<0pt,-7.5pt> ^-{\phi}
\ar@2 "s2"!<0pt,-7.5pt>;"t2"!<0pt,7.5pt> ^-{\psi}
}
\]
Here, we are interested in $3$-polygraphs exclusively, for which we have given a complete classification of critical branchings, see~\cite{GuiraudMalbos09}. They are organised in three families, covering eight different topological configurations of the overlapping, that we will encounter here in different examples. In the case of the $3$-polygraph $\poly{As}_3$, as we have seen in the introduction, we have exactly one critical branching, whose source is an overlapping of two copies of the source~\twocell{(mu *0 1) *1 mu} of the single $3$-cell~\twocell{alpha}:
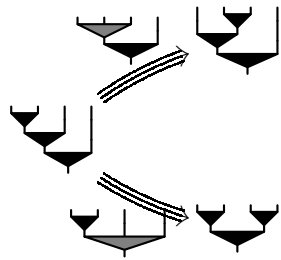
\[
\xymatrix@R=0.5em@C=3em{
& {\twocell{(1 *0 mu *0 1) *1 (mu *0 1) *1 mu}}
\\
{\twocell{(mu *0 2) *1 (mu *0 1) *1 mu}}
	\ar@3@/^/[ur] ^-{\twocell{(alpha *0 1) *1 mu}}
	\ar@3@/_/[dr] _-{\twocell{(mu *0 2) *1 alpha}}
\\
& {\twocell{(mu *0 mu) *1 mu}}
}
\]
The critical branchings are essential in the study of convergence because, under the hypothesis of termination, their confluence ensures the confluence of every branching. This results relies on the fundamental theorem of rewriting theory, namely Newman's lemma, see~\cite{Newman42}, and on another result that depends on the type of rewriting system we consider. The case of $n$-polygraphs is examined in~\cite{GuiraudMalbos09}.

Also, critical branchings of convergent $n$-polygraphs give an algorithmic way to build homotopy bases of track $n$-categories. Indeed, for a given convergent $n$-polygraph $\Sigma$, we define a \emph{basis of generating confluences of $\Sigma$} as a cellular extension of the free $n$-category $\Sigma^*$ made of one $(n+1)$-cell 
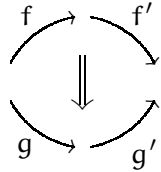
\[
\xymatrix@R=0.5em{
& \strut
	\ar@/^/[dr] ^-{f'} 
	\ar@2 []!<0pt,-15pt>;[dd]!<0pt,15pt>
\\
\strut 
	\ar@/^/[ur] ^-{f}
	\ar@/_/[dr] _-{g}
&& \strut
\\
& \strut
	\ar@/_/[ur] _-{g'}
}
\]
for each critical branching $(f,g)$ of $\Sigma$, where $f'$ and $g'$ are arbitrarily chosen $n$-cells of $\Sigma^*$ with the same target and such that $(f,f')$ and $(g,g')$ are composable. Then we have the following result:

\begin{theorem}[\cite{GuiraudMalbos09}]
\label{TheoremHomotopyBasis}
Let $\Sigma$ be a convergent $n$-polygraph. Then every basis of generating confluences of $\Sigma$ is a homotopy basis of the track $n$-category $\tck{\Sigma}$.
\end{theorem}

\begin{example}
\label{polygrapheAss}
The $3$-polygraph $\poly{As}_3$, seen in the introduction, has one $0$-cell, one $1$-cell, one $2$-cell \twocell{mu} and one $3$-cell
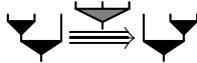
\[
\xymatrix{
{\twocell{ (mu *0 1) *1 mu }}
	\ar@3[r] ^-{\twocell{alpha}}
& {\twocell{ (1 *0 mu) *1 mu } }
}
\]
This $3$-polygraph terminates, see~\cite{GuiraudMalbos09} or the proof of Proposition~\ref{Mon_4HomotopyBase}. It has exactly one critical branching, which is confluent:
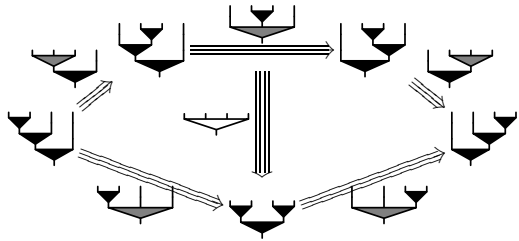
\[
\scalebox{0.8}{
\xymatrix@R=10pt@C=15pt{
&{ \twocell{ (1 *0 mu *0 1) *1 (mu *0 1) *1 mu } }
	\ar@3 [rr] 
		^-{\twocell{(1 *0 mu *0 1) *1 alpha}} 
		_-{}="1" 
&& { \twocell{ (1 *0 mu *0 1) *1 (1 *0 mu) *1 mu } }
	\ar@3 [dr] ^-{\twocell{(1 *0 alpha) *1 mu} }
\\
{ \twocell{ (mu *0 2) *1 (mu *0 1) *1 mu } }
	\ar@3 [ur] ^-{\twocell{(alpha *0 1) *1 mu}}
	\ar@3 [drr] _-{\twocell{(mu *0 2)*1 alpha}}
&&&& { \twocell{ (2 *0 mu) *1 (1 *0 mu) *1 mu } }
\\
&& { \twocell{ (mu *0 mu) *1 mu } }
	\ar@3 [urr] _-{\twocell{(2 *0 mu)*1 alpha}}
	\ar@4 "1"!<0pt,-10pt>;[]!<0pt,20pt> _-{\twocell{aleph}\;}
}
}
\]
As a consequence, filling this $3$-sphere with the $4$-cell \twocell{aleph}, as above, yields a homotopy basis $\poly{As}_4$ of the track $3$-category $\tck{\poly{As}}_3$. In other terms, any two parallel $3$-cells $f$ and $g$ of $\tck{\poly{As}}_3$ are identified in the quotient track $3$-category $\tck{\poly{As}}_3/\poly{As}_4$.
\end{example}

\subsection{Higher-dimensional pro(p)s}
\label{Subsection:HigherDimensionalPRO(P)S}

\subsubsection{Higher-dimensional monoids}

For $n\geq 1$, a \emph{(track) $n$-monoid} is a (track) $n$-category with exactly one $0$-cell, see~\cite{Burroni93}. In particular, a $1$-monoid is a monoid, a track $1$-monoid is a group, a $2$-monoid is a strict monoidal category, a track $2$-monoid is a strict monoidal groupoid. More generally, for $n\geq 2$, a (track) $n$-monoid is a strict monoidal category enriched in (track) $(n-2)$-categories. When the corresponding (enriched) monoidal category is symmetric, we say that an $n$-monoid is \emph{symmetric}.

\subsubsection{Higher-dimensional pro(p)s}

For $n\geq 1$, a \emph{(track) $n$-pro} is a (track) $n$-monoid whose underlying ($1$-)monoid is the monoid $\Nb$ of natural numbers with the addition. A \emph{(track) $n$-prop} is a symmetric (track) $n$-pro. In particular, $2$-pro(p)s coincide with Mac Lane's PRO(P)s, see~\cite{MacLane65}. Let us note that we could consider a more general definition of $n$-pro(p)s by replacing the monoid $\Nb$, which is the free monoid on one generator, by any free monoid.

Here, we are interested in track $3$-pro(p)s and, more precisely, in four main examples: the track $3$-pros $\As\Cat$ of categories with an associative product (see the introduction) and $\Mon\Cat$ of monoidal categories (see~\ref{2ProMonCat}) and the track $3$-props $\Sym\Cat$ of symmetric monoidal categories (see~\ref{2PropSymCat}) and $\Brd\Cat$ of braided monoidal categories (see~\ref{2PropBrCat}). We will show, for each one of those $3$-pro(p)s, how to use homotopy bases built from convergent presentations in order to prove a coherence theorem for the corresponding algebras, a notion we introduce now.

\subsubsection{Algebras over \pdf{3}-pro(p)s}
\label{algebrasOverPro(p)}

We see the (large) monoidal $2$-category $\Cat$ of (small) categories, functors and natural transformations as a (large) $3$-monoid with categories as $1$-cells, functors as $2$-cells, natural transformations as $3$-cells, cartesian product as $0$-composition, composition of functors as $1$-composition, vertical composition of natural transformations as $2$-composition.

If $\P$ is a $3$-pro (resp. $3$-prop), a \emph{$\P$-algebra} is a $3$-functor from $\P$ to $\Cat$ (resp. whose corresponding strict monoidal $2$-functor preserves the symmetry). If $\Ar$ and $\Br$ are $\P$-algebras, a \emph{morphism of $\P$-algebras from $\Ar$ to $\Br$} is a natural transformation from $\Ar$ to $\Br$, \ie, a pair $(F,\Phi)$ where $F:\Ar(1) \fl \Br(1)$ is a functor and $\Phi$ is a map sending every $2$-cell $f:m\dfl n$ in $\P$ to a natural isomorphism with the following shape, where $\Cr=\Ar(1)$ and $\Dr=\Br(1)$,
\[
\xymatrix@R=1em{
& {\Dr^m}
	\ar@/^/ [dr] ^-{\Br(f)}
	\ar@2 []!<0pt,-10pt>;[dd]!<0pt,10pt> _-{\Phi_f} 
\\
{\Cr^m} 
	\ar@/^/ [ur] ^-{F^m}
	\ar@/_/ [dr] _-{\Ar(f)}
&& {\Dr^n}
\\
& {\Cr^n}
	\ar@/_/ [ur] _-{F^n}
}
\]
such that the following relations hold:
\begin{itemize}
\item for every $2$-cells $f:m\dfl n$ and $g:p\dfl q$ of $\P$, we have $\Phi_{f\star_0 g} = \Phi_f\times \Phi_g$:
\[
\xymatrix@R=1em{
& {\Dr^{m+p}}
	\ar@/^/ [dr] ^-{\Br(f\star_0 g)}
	\ar@2 []!<0pt,-10pt>;[dd]!<0pt,10pt> |-{\Phi_{f\star_0 g}} 
&&&& {\Dr^m\times\Dr^p}
	\ar@/^/ [dr] ^-{\Br(f)\times\Br(g)}
	\ar@2 []!<0pt,-10pt>;[dd]!<0pt,10pt> |-{\Phi_{f}\times\Phi_{g}} 
\\
{\Cr^{m+p}} 
	\ar@/^/ [ur] ^-{F^{m+p}}
	\ar@/_/ [dr] _-{\Ar(f\star_0 g)}
&& {\Dr^{n+q}}
& = 
& {\Cr^m\times\Cr^p}
	\ar@/^/ [ur] ^-{F^{m}\times F^{p}}
	\ar@/_/ [dr] _-{\Ar(f)\times\Ar(g)}
&& {\Dr^{n}\times\Dr^{q}}
\\
& {\Cr^{n+q}}
	\ar@/_/ [ur] _-{F^{n+q}}
&&&& {\Cr^{n}\times\Cr^{q}}
	\ar@/_/ [ur] _-{F^{n}\times F^{q}}
}
\]
\item for every $2$-cells $f:m\dfl n$ and $g:n\dfl p$ in $\P$, we have $\Phi_{f\star_1 g} = (\Phi_f\star_1 \Br(g)) \star_2 (\Ar(f)\star_1 \Phi_g)$:
\[
\xymatrix@R=1em{
& {\Dr^{m}}
	\ar@/^/ [dr] ^-{\Br(f\star_1 g)}
	\ar@2 []!<0pt,-10pt>;[dd]!<0pt,10pt> |-{\Phi_{f\star_1 g}}  
\\
{\Cr^{m}} 
	\ar@/^/ [ur] ^-{F^{m}}
	\ar@/_/ [dr] _-{\Ar(f\star_1 g)}
&& {\Dr^{p}}
\\
& {\Cr^{p}}
	\ar@/_/ [ur] _-{F^{p}}
}
\qquad \raisebox{-5ex}{=} \qquad
\raisebox{5ex}{
\xymatrix@R=1em{
& {\Dr^m}
	\ar@/^/ [dr] ^-{\Br(f)}
	\ar@2 []!<0pt,-10pt>;[dd]!<0pt,10pt> _-{\Phi_{f}}  
\\
{\Cr^m} 
	\ar@/^/ [ur] ^-{F^{m}}
	\ar@/_/ [dr] _-{\Ar(f)}
&& {\Dr^n}
	\ar@/^/ [dr] ^-{\Br(g)}
	\ar@2 []!<0pt,-10pt>;[dd]!<0pt,10pt> ^-{\Phi_{g}}  
\\ 
& {\Cr^n}
	\ar@/_/ [dr] _-{\Ar(g)}
	\ar [ur] |-{F^n}
&& {\Dr^p}
\\
&& {\Cr^p}
	\ar@/_/ [ur] _-{F^p}
}
}
\]
\item for every $3$-cell $A:f\tfl g:m\dfl n$ in $\P$, we have $\Phi_f\star_2 (\Ar(A) \star_1 F^n) \:=\: (F^m \star_1 \Br(A)) \star_2 \Phi_g$:
\[
\xymatrix@R=1em{
& {\Dr^m}
	\ar@/^/ [dr] ^-{\Br(f)}
	\ar@2 []!<0pt,-10pt>;[dd]!<0pt,10pt> ^-{\Phi_f}  
\\
{\Cr^m} 
	\ar@/^/ [ur] ^-{F^m}
	\ar [dr] ^(0.4){\Ar(f)} _-{}="src"
	\ar@/_7ex/ [dr] _-{\Ar(g)} ^-{}="tgt"
	\ar@2 "src";"tgt" |-{\Ar(A)}
&& {\Dr^n}
\\
& {\Cr^n}
	\ar@/_/ [ur] _-{F^n}
}
\qquad \raisebox{-5ex}{=} \qquad
\xymatrix@R=1em{
& {\Dr^m}
	\ar@/^7ex/ [dr] ^-{\Br(f)} _-{}="src"
	\ar [dr] _(0.6){\Br(g)} ^-{}="tgt"
	\ar@2 []!<0pt,-10pt>;[dd]!<0pt,10pt> _-{\Phi_g} 
	\ar@2 "src";"tgt" |-{\Br(A)}
\\
{\Cr^m} 
	\ar@/^/ [ur] ^-{F^m}
	\ar@/_/ [dr] _-{\Ar(g)}
&& {\Dr^n}
\\
& {\Cr^n}
	\ar@/_/ [ur] _-{F^n}
}
\]
\end{itemize}
The $\P$-algebras and their morphisms form a category, denoted by $\Alg(\P)$. 

\subsubsection{Coherence problem for algebras over a \pdf{3}-pro(p)}
\label{coherenceProblemPro(p)}

Let $\P$ be a $3$-pro(p) and let $\Ar$ be a $\P$-algebra. A \emph{$\P$-diagram in $\Ar$} is the image~$\Ar(\gamma)$ of a $3$-sphere $\gamma$ in $\P$. A $\P$-diagram $\Ar(\gamma)$ in $\Ar$ \emph{commutes} if the relation $\Ar(s(\gamma))= \Ar(t(\gamma))$ is satisfied in $\Cat$. The coherence problem for algebras over a $3$-pro(p) is:

\medskip
\noindent
\textsc{Coherence problem:} 

\begin{minipage}{14cm}
\emph{Given a $3$-pro(p) $\P$, does every $\P$-diagram commute in every $\P$-algebra?}
\end{minipage}

\medskip
\noindent
As a consequence of the definition of an aspherical $3$-pro(p), we have the following sufficient condition for giving a positive answer to the coherence problem:

\begin{proposition}
\label{Propaspherical=>coherence}
If $\P$ is an aspherical $3$-pro(p) then every $\P$-diagram commutes in every $\P$-algebra.
\end{proposition}

\begin{example}
Let $\As\Cat$ be the track $3$-pro defined as the following quotient:
\[
\As 
	\:=\: \tck{\poly{As}}_3 \;\big/\; \poly{As}_4 
	\:=\: \tck{\left(\:\twocell{1}\:, \twocell{mu}, \twocell{alpha}\right)} \;\big/\; \twocell{aleph}.
\]
The category $\Alg(\As\Cat)$ is isomorphic to the category of (small) associative categories, the correspondence between an associative category $(\Cr,\tens,\alpha)$ and a $3$-functor $\Ar:\As\Cat\fl\Cat$ being given by
\[
\Ar\left(\:\twocell{1}\:\right) \:=\: \Cr, 
\qquad\qquad 
\Ar\left(\twocell{mu}\right) \:=\: \tens,
\qquad\qquad 
\Ar\left(\twocell{alpha}\right) \:=\: \alpha.
\]
This correspondence is well-defined since the coherence diagram satisfied by associative categories corresponds to the $4$-cell~\twocell{aleph}. We have seen that $\poly{As}_4=\ens{\twocell{aleph}}$ is a homotopy basis of $\tck{\poly{As}}_3$, so that $\As\Cat$ is an aspherical track $3$-pro. As a consequence, in every associative category $\Cr$, every $\As\Cat$-diagram is commutative. This fact can be informally restated as: every diagram built in $\Cr$ from the functor $\tens$ and the natural transformation~$\alpha$ is commutative. 
\end{example}

\section{Coherence in monoidal categories}
\label{Section:CoherenceMonoidalCategories}

\subsection{Coherence in algebras over track \pdf{3}-pros}

\subsubsection{Presentations of track \pdf{3}-pros}

Let $\P$ be a track $3$-pro. A \emph{presentation of $\P$} is a pair $(\Sigma_3,\Sigma_4)$, where~$\Sigma_3$ is a $3$-polygraph and~$\Sigma_4$ is a cellular extension of the free track $3$-category $\tck{\Sigma}_3$ such that 
\[
\P \:\simeq\: \tck{\Sigma}_3/\Sigma_4.
\]
Note that, in that case, the $3$-polygraph $\Sigma_3$ has exactly one $0$-cell and one $1$-cell. 

A presentation of a track $3$-pro yields a diagram which is similar to the one corresponding to the inductive construction of a $4$-polygraph, see Section~\ref{Polygraphs}: 
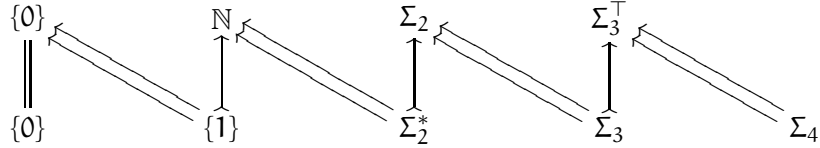
\[
\xymatrix{
\ens{0} 
     \ar@{=}[d] 
&& \Nb 
&& \Sigma_2
&& \tck{\Sigma}_3
\\
\ens{0}
&& \ens{1}
     \ar@<0.5ex>[llu]
     \ar@<-0.5ex>[llu]
     \ar@{>->}[u]
&& \Sigma_2^*
     \ar@<0.5ex>[llu]
     \ar@<-0.5ex>[llu]
     \ar@{>->}[u]
&& \Sigma_3
     \ar@<0.5ex>[llu]
     \ar@<-0.5ex>[llu]
     \ar@{>->}[u]
&& \Sigma_4
     \ar@<0.5ex>[llu]
     \ar@<-0.5ex>[llu]
}
\]
A presentation $(\Sigma_3,\Sigma_4)$ of $\P$ is \emph{convergent} when $\Sigma_3$ is a convergent $3$-polygraph and $\Sigma_4$ is a cellular extension of generating confluences of~$\Sigma_3$. 

By definition, $\P$ is an aspherical $3$-pro if and only if, for every presentation $(\Sigma_3,\Sigma_4)$ of $\P$, the cellular extension $\Sigma_4$ is a homotopy basis. The latter condition is satisfied by any convergent presentation of~$\P$, yielding the following sufficient condition for giving a positive answer to the coherence problem for $\P$-algebras:

\begin{theorem}
\label{TheoremProAspherical}
If a track $3$-pro $\P$ admits a convergent presentation then every $\P$-diagram commutes in every $\P$-algebra.
\end{theorem}

\subsection{Identities among relations for presentations of track \pdf{3}-pros}

This section is based on notions and results from~\cite{GuiraudMalbos10smf}, that we briefly recall first. 

\subsubsection{Contexts and natural systems} 

Let $\C$ be an $n$-category. A \emph{context of $\C$} is an $(n+1)$-cell $C$ of some free $(n+1)$-category $\C[x]$, where $x$ is an $n$-sphere of $\C$, such that $C$ contains exactly one occurrence of $x$. If $f$ is an $n$-cell of $\C$ which is parallel to $x$, we denote by $C[f]$ the $n$-cell of $\C$ obtained by replacing~$x$ with $f$ in $C$ and, if $D$ is a context which is parallel to $x$, we denote by $C\circ D$ the context of $\C$ obtained by replacing $x$ with $D$ in $C$. 

A \emph{whisker of $\C$} is a context of $\C$ that contains only $(n-1)$-cells of $\C$, apart from the $n$-sphere $x$. Note that whiskers of $\C$ are in bijective correspondence with contexts of the $(n-1)$-category $\C_{n-1}$ underlying~$\C$.

The contexts of $\C$ form a category whose objects are the $n$-cells of $\C$ and whose morphisms from~$f$ to~$g$ are the contexts $C$ of $\C$ such that $C[f]=g$. A \emph{natural system on $\C$} is a functor from the category of contexts of $\C$ to the category of abelian groups.

\subsubsection{Abelian track \pdf{n}-category}

Let $\T$ be a track $n$-category. An $n$-cell $f$ of $\T$ is \emph{closed} when its source and its target are equal; this common $(n-1)$-cell is the \emph{base cell of $f$}. For every $(n-1)$-cell $u$ of $\T$, the $n$-cells of $\T$, equipped with the composition $\star_{n-1}$, form a group, which is denoted by $\Aut^{\T}_u$.

We say that a track $n$-category $\T$ is \emph{abelian} when every group $\Aut^{\T}_u$ is abelian, \ie, when, for every closed $n$-cells $f$ and $g$ with same base cell, the relation $f\star_{n-1}g=g\star_{n-1} f$ is satisfied. Note that, for an abelian track $n$-category $\T$, the assignment of each $(n-1)$-cell $u$ of $\T$ to the abelian group $\Aut^{\T}_u$ extends to a natural system $\Aut^{\T}$ on the $(n-1)$-category $\T_{n-1}$ underlying $\T$.

We denote by $\ab{\T}$ the \emph{abelianised track $n$-category of $\T$}, defined as the quotient of $\T$ by the cellular extension made of one $(n+1)$-cell from $f\star_{n-1} g$ to $g\star_{n-1} f$ for every pair $(f,g)$ of closed $n$-cells of $\T$ with the same base cell.

\subsubsection{Identities among relations for \pdf{n}-polygraphs} 

Let $\Sigma$ be an $n$-polygraph. We denote by $\cl{u}$ the image of an $(n-1)$-cell $u$ of $\abtck{\Sigma}$ by the canonical projection to the $(n-1)$-category $\cl{\Sigma}$ presented by $\Sigma$. We define the natural system on $\cl{\Sigma}$ of \emph{identities among relations of $\Sigma$}, denoted by $\Pi(\Sigma)$, as follows. 

For any $(n-1)$-cell $u$ in $\cl{\Sigma}$, the abelian group $\Pi(\Sigma)_{u}$ is defined as the group with one generator $\iar{f}$ for every $n$-cell $f:v\fl v$ of $\abtck{\Sigma}$ with $\cl{v}=u$, subjected to the following relations:
\begin{enumerate}[\bf i)]
\item $\iar{f\star_{n-1} g} = \iar{f} + \iar{g}$, for every $n$-cells $f,g:v\fl v$ of $\abtck{\Sigma}$ with $\cl{v}=u$;  
\item $\iar{f\star_{n-1} g} = \iar{g\star_{n-1}f}$, for every $n$-cells $f:v\fl w$ and $g:w\fl v$ of $\abtck{\Sigma}$ with $\cl{v}=\cl{w}=u$.  
\end{enumerate}

\noindent
For any context $C$ of $\cl{\Sigma}$ from $u$ to $v$, the morphism of groups $\Pi(\Sigma)_C$ from $\Pi(\Sigma)_u$ to $\Pi(\Sigma)_v$ is given, on a generator $\iar{f}$, by $\Pi(\Sigma)_C (\iar{f}) = \iar{\rep{C}[f]}$, where $\rep{C}$ is any whisker of $\abtck{\Sigma}$ that represents the context $C$ of~$\cl{\Sigma}$. 

In~\cite{GuiraudMalbos10smf}, the authors prove that the natural system $\Pi(\Sigma)$ is well-defined and, in particular, that its values on contexts do not depend on the chosen representatives. Moreover, the functor $\Pi(\Sigma)$ is the unique natural system on $\cl{\Sigma}$, up to isomorphism, such that there exists an isomorphism of natural systems on~$\Sigma^*_{n-1}$
\[
\Phi \::\: \rep{\Pi(\Sigma)} \longrightarrow \Aut^{\abtck{\Sigma}},
\]
where $\rep{\Pi(\Sigma)}$ is defined, on an $(n-1)$-cell $u$ of $\abtck{\Sigma}$, by $\rep{\Pi(\Sigma)}_u \:=\: \Pi(\Sigma)_{\cl{u}}$. The isomorphism $\Phi$ is given, for an $(n-1)$-cell $u$ of $\abtck{\Sigma}$ and a closed $n$-cell $f$ of $\Pi(\Sigma)$ with base $v$ such that $\cl{u}=\cl{v}$, by
\[
\Phi(\iar{f}) \:=\: g \star_{n-1} f \star_{n-1} g^-,
\]
where $g:u\fl v$ is any $n$-cell of $\abtck{\Sigma}$. 

Let $\Gamma$ be a cellular extension of $\tck{\Sigma}$. For each $\gamma$ in $\Gamma$, we denote by $\tilde{\gamma}$ the following $n$-cell of $\tck{\Sigma}$:
\[
\tilde{\gamma} \:=\: s(\gamma)\star_{n-1} t(\gamma)^-.
\]
We define $\tilde{\Gamma}=\ens{\tilde{\gamma},\gamma\in\Gamma}$. When $\Gamma$ is a homotopy basis, then the set  $\iar{\tilde{\Gamma}}$ is a generating set of the natural system~$\Pi(\Sigma)$, \ie, every element $a$ of any $\Pi(\Sigma)_u$ can be written
\[
a \:=\: \sum_{i=1}^n \epsilon_i C_i\iar{b_i},
\]
where each $b_i$ is a element of $\tilde{\Gamma}$, each $C_i$ is a context of $\cl{\Sigma}$ and each $\epsilon_i$ is $\pm 1$. The proof relies on an equivalence between the facts that $\Gamma$ is a homotopy basis and that every $n$-cell $f$ of $\tck{\Sigma}$ can be written
\[
f \:=\: 
	\left(g_1\star_{n-1} C_1 [\tilde{\gamma}_1^{\epsilon_1}] \star_{n-1} g_1^-\right) 
	\star_{n-1} \cdots\star_{n-1}
	\left(g_k\star_{n-1} C_k [\tilde{\gamma}_k^{\epsilon_k}] \star_{n-1} g_k^-\right), 
\]
where each $\gamma_i$ is in $\Gamma$, each $\epsilon_i$ is $\pm 1$, each $C_i$ is a whisker of $\tck{\Sigma}$ and each $g_i$ is an $n$-cell of $\tck{\Sigma}$. We refer the reader to~\cite{GuiraudMalbos10smf} for the proof. Here, in the special case of presentations of track $3$-pros, we get:

\begin{proposition}
\label{generatingSetIar}
If $\P$ is an aspherical track $3$-pro then, for every presentation $(\Sigma_3,\Sigma_4)$ of $\P$, the natural system $\Pi(\Sigma_3)$ on the $2$-pro $\cl{\Sigma}_3$ is generated by the set~$\iar{\tilde{\Sigma}_4}$.
\end{proposition}

\begin{proof}
Since $\P$ is aspherical, then $\Sigma_4$ is a homotopy basis of the track $3$-category $\tck{\Sigma}_3$ and, thus, of the abelianised track $3$-category $\ab{(\tck{\Sigma}_3)}$. Hence, any closed $3$-cell~$A$ in $\ab{(\tck{\Sigma}_3)}$ can be written  
\[
A \:=\: \big( A_1\star_2 C_1[\tilde{\omega}_1^{\epsilon_1}] \star_2 A_1^- \big) 
	\star_2 \dots \star_2
	\big( A_k\star_2 C_k[\tilde{\omega}_k^{\epsilon_k}] \star_2 g_k^- \big),
\]
where each $\omega_i$ is in $\Sigma_4$, each $\epsilon_i$ is $\pm 1$, each $C_i$ is a whisker of $\ab{(\tck{\Sigma}_3)}$ and each $A_i$ is a $3$-cell of $\ab{(\tck{\Sigma}_3)}$. 

Hence, any generator $\iar{A}$ of $\Pi(\Sigma_3)_f$, for $f$ a $2$-cell of $\cl{\Sigma}_3$, can be written
\[
\iar{A} 
	\:=\: \sum_{i=1}^k \iar{A_i \star_2 C_i[\tilde{\omega}_i^{\epsilon_i}] \star_2 A_i^-} 
	\:=\: \sum_{i=1}^k \epsilon_i C_i\iar{\tilde{\omega}_i}.
\]
Thus, the elements of $\iar{\tilde{\Sigma}_4}$ form a generating set for $\Pi(\Sigma_3)$. 
\end{proof}

\begin{corollary}
Let $\P$ be an aspherical track $3$-pro and let $(\Sigma_3,\Sigma_4)$ be a presentation of $\P$. If $\Sigma_4$ is finite, then the natural system $\Pi(\Sigma_3)$ is finitely generated.
\end{corollary}

\subsection{Application: coherence for monoidal categories}
\label{Subsection:ApplicationCoherenceMon}

We recall that a \emph{monoidal category} is a category $\Cr$, equipped with two functors $\tens:\Cr\times\Cr\fl\Cr$ and $e:\ast\fl\Cr$, and three natural isomorphisms
\[
\alpha_{x,y,z} \: : \: (x\tens y)\tens z \:\fl\: x\tens(y\tens z) \:,
\qquad
\lambda_x \::\: e \tens x \:\fl\: x,
\qquad
\rho_x \::\: x\tens e \:\fl\: x,
\]
such that the following two diagrams commute in $\Cr$:
\[
\scalebox{0.9}{
\xymatrix@C=-3em{
& { (x\tens (y\tens z))\tens t}
	\ar[rr] ^-{\alpha}
&& { x\tens ((y\tens z) \tens t)}
	\ar[dr] ^-{\alpha}
\\
{ ((x\tens y)\tens z) \tens t}
	\ar[ur] ^-{\alpha}
	\ar[drr] _-{\alpha}
&& { \copyright}
&& { x\tens(y\tens (z\tens t))}
\\
&& { (x\tens y)\tens (z\tens t)}
	\ar[urr] _-{\alpha}
}
\raisebox{-2em}{
\xymatrix@C=0.5em{
& { x\tens (e\tens y)}
	\ar[dr] ^-{\lambda}
\\
{ (x\tens e)\tens y}
	\ar[ur] ^-{\alpha}
	\ar[rr] _-{\rho} ^-{}="1"
&& { x\tens y}
\ar@{} "1,2" ; "1" |-{ \copyright}
}
}
}
\]
A \emph{monoidal functor} from $\Cr$ to $\Dr$ is a triple $(F,\upphi,\iota)$ made of a functor $F:\Cr\fl\Dr$ and two natural natural isomorphisms $\upphi_{x,y} : Fx\tens Fy \fl F(x\tens y)$ and $\iota : e \fl F(e)$ such that the following diagrams commute in $\Dr$:
\[
\label{relationFunctor1}
\scalebox{0.9}{
\xymatrix@C=1em{
& Fx\tens (Fy\tens Fz) 
	\ar[rr] ^-{1\tens\upphi}
&& Fx\tens F(y \tens z)
	\ar[dr] ^-{\upphi}
\\
(Fx \tens Fy) \tens Fz
	\ar[ur] ^-{\alpha}
	\ar[dr] _-{\upphi\tens 1}
	\ar@{} [rrrr] |-{\copyright}
&&&& F(x\tens(y\tens z))
\\
& F(x\tens y)\tens Fz 
	\ar[rr] _-{\upphi}
&& F((x\tens y)\tens z)
	\ar[ur] _-{F\alpha}
}
}
\]
\[
\label{relationFunctor2}
\scalebox{0.9}{
\xymatrix{
{ Fx\tens e }
	\ar[r] ^-{\rho}
        \ar[d] _-{1\tens\iota}
        \ar@{} [dr] |-{\copyright}
&
{ Fx}
\\
{ Fx\tens Fe}
	\ar[r] _-{\upphi}
&
{ F(x\tens e)}
        \ar[u] _-{F\rho}
}
\qquad \qquad
\xymatrix{
{ e\tens Fx}
	\ar[r] ^-{\lambda}
        \ar[d] _-{\iota\tens 1}
        \ar@{} [dr] |-{\copyright}
&
{ Fx}
\\
{ Fe\tens Fx}
	\ar[r] _-{\upphi}
&
{ F(e\tens x)}
        \ar[u] _-{F\lambda}
}
}
\]

\subsubsection{The \pdf{3}-pro of monoidal categories}
\label{2ProMonCat}

Let $\Mon\Cat$ be the $3$-pro presented by $(\poly{Mon}_3,\poly{Mon}_4)$, where $\poly{Mon}_3$ is the $3$-polygraph with two $2$-cells~\twocell{mu},~\twocell{eta} and three $3$-cells 
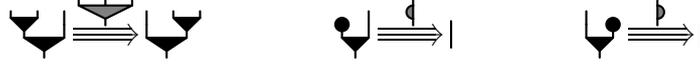
\[
\label{Mon_3}
\xymatrix{
{\twocell{ (mu *0 1) *1 mu }}
	\ar@3[r] ^-{\twocell{alpha}}
& {\twocell{ (1 *0 mu) *1 mu } }
}
\qquad\qquad
\xymatrix{
{\twocell{ (eta *0 1) *1 mu }}
	\ar@3[r] ^-{\twocell{lambda}}
& {\twocell{1}}
}  
\qquad\qquad
\xymatrix{
{\twocell{ (1 *0 eta) *1 mu }}
	\ar@3[r] ^-{\twocell{rho}}
& {\twocell{1}}
}
\]
and $\poly{Mon}_4$ is the cellular extension of $\tck{\poly{Mon}}_3$ made of the following two $4$-cells: 
\[
\label{Mon_4}
\scalebox{0.8}{
\xymatrix@R=10pt@C=15pt{
&{ \twocell{ (1 *0 mu *0 1) *1 (mu *0 1) *1 mu } }
	\ar@3 [rr] 
		^-{\twocell{(1 *0 mu *0 1) *1 alpha}} 
		_-{}="1" 
&& { \twocell{ (1 *0 mu *0 1) *1 (1 *0 mu) *1 mu } }
	\ar@3 [dr] ^-{\twocell{(1 *0 alpha) *1 mu} }
\\
{ \twocell{ (mu *0 2) *1 (mu *0 1) *1 mu } }
	\ar@3 [ur] ^-{\twocell{(alpha *0 1) *1 mu}}
	\ar@3 [drr] _-{\twocell{(mu *0 2)*1 alpha}}
&&&& { \twocell{ (2 *0 mu) *1 (1 *0 mu) *1 mu } }
\\
&& { \twocell{ (mu *0 mu) *1 mu } }
	\ar@3 [urr] _-{\twocell{(2 *0 mu)*1 alpha}}
	\ar@4 "1"!<0pt,-10pt>;[]!<0pt,20pt> _-*+{\twocell{aleph}}
}
}
\qquad\qquad
\scalebox{0.8}{
\xymatrix{
& { \twocell{ ( 1 *0 eta *0 1) *1 (1 *0 mu) *1 mu } }
	\ar@3 [dr] ^-{\twocell{(1 *0 lambda) *1 mu}}
\\
{ \twocell{ ( 1 *0 eta *0 1) *1 (mu *0 1) *1 mu } }
	\ar@3 [ur] ^-{\twocell{(1 *0 eta *0 1) *1 alpha}}
	\ar@3 [rr] _-{\twocell{(rho *0 1) *1 mu}} ^-{}="1"
&& { \twocell{mu} }
\ar@4 "1,2"!<0pt,-20pt>;"1"!<0pt,10pt> _-*+{\twocell{beth}}
}
}
\]

\begin{lemma}
\label{isoAlgMon}
The category of small monoidal categories and monoidal functors is isomorphic to the category $\Alg(\Mon\Cat)$.
\end{lemma}

\begin{proof}
For a monoidal category $(\Cr,\tens,e,\alpha,\lambda,\rho)$, the corresponding $\Mon\Cat$-algebra $\Ar$ is given by:
\begin{equation}
\label{isoAlgMonEq}
\Ar(\:\twocell{1}\:) \:=\: \Cr,
\quad
\Ar(\twocell{mu}) \:=\: \tens,
\quad
\Ar(\twocell{eta}) \:=\: e,
\quad
\Ar(\twocell{alpha}) \:=\: \alpha,
\quad
\Ar(\twocell{lambda}) \:=\: \lambda,
\quad
\Ar(\twocell{rho}) \:=\: \rho.
\end{equation}
The two commutative diagrams satisfied by monoidal categories correspond to the $\Mon\Cat$-diagrams $\Ar(\twocell{aleph})$ and $\Ar(\twocell{beth})$. If~$(F,\upphi,\iota)$ is a monoidal functor, the corresponding morphism~$\Psi$ of $\Mon\Cat$-algebras is: 
\[
\Psi_{\twocell{1}} \:=\: F,
\qquad
\Psi_{\twocell{mu}} \:=\: \upphi,
\qquad
\Psi_{\twocell{eta}} \:=\: \iota.
\qedhere
\]  
\end{proof}

\begin{proposition}[\cite{GuiraudMalbos09}]
\label{Mon_4HomotopyBase}
The cellular extension $\poly{Mon}_4$ of the free track $3$-category $\tck{\poly{Mon}}_3$ is a homotopy basis. 
\end{proposition}

\begin{proof}
First, we check that the $3$-polygraph $\poly{Mon}_3$ terminates. We recall the proof from~\cite{Guiraud06jpaa}, see also~\cite{GuiraudMalbos09}. We consider the $2$-functor $X$ from $\poly{Mon}^*_2$ to the category of ordered sets and monotone maps, seen as a $2$-category with one $0$-cell:
\[
X(\:\twocell{1}\:) \:=\: \Nb\setminus\ens{0}, \qquad
X(\twocell{mu}) (i,j) \:=\: i+j, \qquad
X(\twocell{eta}) \:=\: 1.
\]
Then, we consider the following assignment of $2$-cells of $\poly{Mon}_2$:
\[
\dr (\twocell{mu}) (i,j) \:=\: i, \qquad 
\dr (\twocell{eta}) \:=\: 0.
\]
This assignment extends, in a unique way, to a derivation of $\poly{Mon}_2^*$ with values in $X$, \ie, a map $\dr$ that sends each $2$-cell $f:m\dfl n$ of $\poly{Mon}_2^*$ to a monotone map $\dr(f):\Nb^m\fl\Nb$ that satisfies the following relations:
\[
\dr(f\star_0 g)(i_1,\dots,i_{m+n}) \:=\: \dr(f)(i_1,\dots,i_m) \:+\: \dr(g)(i_{m+1},\dots,i_{m+n})
\]
and
\[
\dr(f\star_1 g)(i_1,\dots,i_m) \:=\: \dr(f)(i_1,\dots,i_m) + \dr(g)\circ X(f)(i_1,\dots,i_m).
\]
We check that, for every $3$-cell $\alpha$ of $\poly{Mon}_3$, we have
\[
X(s(\alpha)) \:\geq\: X(t(\alpha)) 
\qquad\text{and}\qquad
\dr(s(\alpha)) \:>\: \dr(t(\alpha)),
\]
where monotone maps are compared pointwise. This implies that, for every non-degenerate $3$-cell $A$ of $\poly{Mon}_3^*$, we have $\dr(s(A))>\dr(t(A))$. Since $\dr$ takes its values in $\Nb$, the $3$-polygraph $\poly{Mon}_3$ terminates.

For confluence, we study the critical branchings of $\poly{Mon}_3$: it has five critical branchings and each of them is confluent. This yields a cellular extension $\Gamma$ of $\tck{\poly{Mon}}_3$ with five $4$-cells, the ones of $\poly{Mon}_4$ plus the following three $4$-cells: 
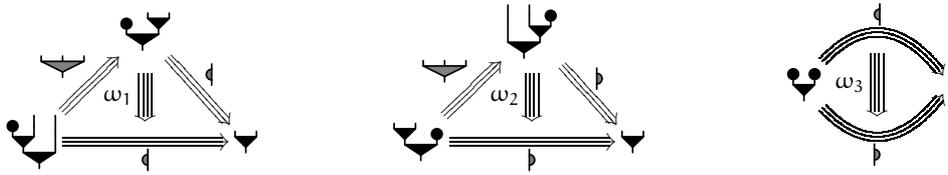
\[
\scalebox{0.8}{
\xymatrix{
& {\twocell{(eta *0 mu) *1 mu}}
	\ar@3 [dr] ^-{\twocell{lambda}}
\\
{\twocell{(((eta *0 1) *1 mu) *0 1) *1 mu}}
	\ar@3 [ur] ^-{\twocell{alpha}}
	\ar@3 [rr] _-{\twocell{lambda}} ^-{}="1"
&& {\twocell{mu}}
\ar@4 "1,2"!<0pt,-20pt>;"1"!<0pt,10pt> _-*+{\omega_1}
}
}
\qquad\qquad
\scalebox{0.8}{
\xymatrix{
& {\twocell{(1 *0 ((1 *0 eta) *1 mu) ) *1 mu }}
	\ar@3 [dr] ^-{\twocell{rho}}
\\
{\twocell{(mu *0 eta) *1 mu}}
	\ar@3 [ur] ^-{\twocell{alpha}}
	\ar@3 [rr] _-{\twocell{rho}} ^-{}="1"
&& {\twocell{mu}}
\ar@4 "1,2"!<0pt,-20pt>;"1"!<0pt,10pt> _-*+{\omega_2}
}
}
\qquad\qquad
\raisebox{-20pt}{
\scalebox{0.8}{
\xymatrix{
{\twocell{(eta *0 eta) *1 mu}}
	\ar@3@/^5ex/ [rr] ^{\twocell{lambda}} _{}="1"
	\ar@3@/_5ex/ [rr] _{\twocell{rho}} ^{}="2"
&& {\twocell{1}}
\ar@4 "1"!<0pt,-10pt>;"2"!<0pt,10pt> _-*+{\omega_3}
}
}
}
\]
Hence $\Gamma$ is a homotopy basis of $\tck{\poly{Mon}}_3$. To prove that $\poly{Mon}_4$ is a homotopy basis, we show that, for each $4$-cell $\omega_i$, we have $\cl{s(\omega_i)}= \cl{t(\omega_i)}$ in $\Mon\Cat$. For $\omega_1$, we define the $4$-cell $\gamma$ of $\tck{\poly{Mon}}_3(\poly{Mon}_4)$ by the following relation, where we abusively denote $3$-cells by the generating $3$-cell of $\poly{\poly{Mon}}_3$ they contain: 
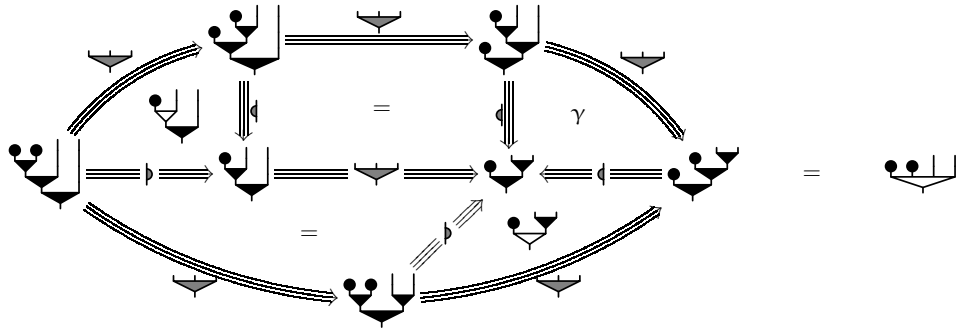
\[
\scalebox{0.8}{
\xymatrix{
&& {\twocell{(eta *0 2) *1 (eta *0 mu *0 1) *1 (mu *0 1) *1 mu}}
	\ar@3 [rr] ^-{\twocell{alpha}} 
	\ar@3 [d] ^-{\twocell{lambda}}
	\ar@{} [drr] |-{=}
&& {\twocell{(eta *0 2) *1 (mu *0 1) *1 (eta *0 mu) *1 mu}}
	\ar@3@/^3ex/ [drr] ^-{\twocell{alpha}} _-{}="2"
	\ar@3 [d] _-{\twocell{lambda}}
\\
 {\twocell{(eta *0 eta *0 2) *1 (mu *0 2) *1 (mu *0 1) *1 mu}}
	\ar@3@/^3ex/ [urr] ^-{\twocell{alpha}} _-{}="1"
	\ar@3 [rr] |-*+[o]{\twocell{rho}} 
	\ar@3@/_3ex/ [drrr] _-{\twocell{alpha}}
&&	{\twocell{(eta *0 2) *1 (mu *0 1) *1 mu}}
	\ar@3 [rr] |-*+[o]{\twocell{alpha}}
	\ar@{} [dr] |-{=}
&& {\twocell{(eta *0 mu) *1 mu}}
&& {\twocell{(eta *0 mu) *1 (eta *0 mu) *1 mu}}
	\ar@3 [ll] |-*+[o]{\twocell{lambda}}
& =
& {\twocell{(eta *0 eta *0 2) *1 aleph}}
\\
&&& {\twocell{(eta *0 eta *0 2) *1 (mu *0 mu) *1 mu}}
	\ar@3 [ur] |-*+[o]{\twocell{rho}}
	\ar@3@/_3ex/ [urrr] _-{\twocell{alpha}} ^-{}="3"
\ar@{} "1" ; "2,3" |-{\twocell{(eta *0 2) *1 (beth *0 1) *1 mu}}
\ar@{} "2" ; "2,5" |-{\gamma}
\ar@{} "3" ; "2,5" |-{\twocell{(eta *0 mu) *1 beth}}
}
}
\]
As a consequence of this construction, we have $\cl{s(\gamma)}=\cl{t(\gamma)}$ in $\Mon\Cat$. Then we build the following diagram, proving that $\cl{s(\omega_1)}=\cl{t(\omega_1)}$ also holds: 
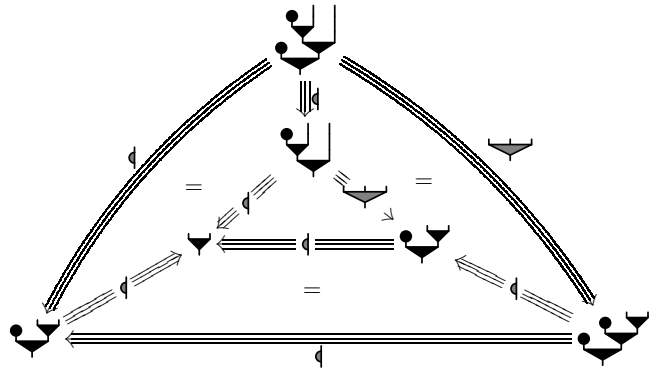
\[
\scalebox{0.8}{
\xymatrix@R=1.5em{
&&& {\twocell{(eta *0 2) *1 (mu *0 1) *1 (eta *0 mu) *1 mu}}
	\ar@3@/_3ex/ [dddlll] _-{\twocell{lambda}} ^-{}="s1"
	\ar@3 [d] ^-{\twocell{lambda}}
	\ar@3@/^3ex/ [dddrrr] ^-{\twocell{alpha}} _-{}="s2"
\\ 
&&& {\twocell{(eta *0 2) *1 (mu *0 1) *1 mu}}
	\ar@3 [dl] |-*+[o]{\twocell{lambda}} _-{}="t1"
	\ar@3 [dr] |-*+[o]{\twocell{alpha}} ^-{}="t2"
\\ 
&& {\twocell{mu}}
&& {\twocell{(eta *0 mu) *1 mu}}
	\ar@3 [ll] |-*+[o]{\twocell{\lambda}} ^-{}="s3"
\\
{\twocell{(eta *0 mu) *1 mu}}
	\ar@3 [urr] |-*+[o]{\twocell{lambda}}
&&&&&& {\twocell{(eta *0 mu) *1 (eta *0 mu) *1 mu}}
	\ar@3 [ull] |-*+[o]{\twocell{lambda}}
	\ar@3 [llllll] ^-{\twocell{lambda}} _-{}="t3"
\ar@{} "s1" ; "t1" |-{=}
\ar@{} "s2" ; "t2" |-{=}
\ar@{} "s3" ; "t3" |-{=}
}
}
\]
For the $4$-cell $\omega_2$, one proceeds in a similar way, starting with the $4$-cell~\twocell{(2 *0 eta *0 eta) *1   aleph}. 

Finally, let us consider the case of the $4$-cell $\omega_3$. We define the $4$-cell $\delta$ of $\tck{\poly{Mon}}_3(\poly{Mon}_4)$ by the following relation: 
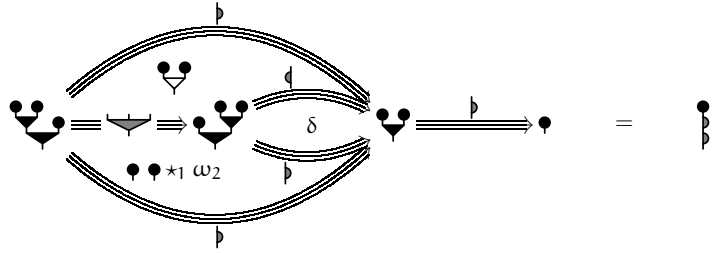
\[
\scalebox{0.8}{
\xymatrix{
{\twocell{(eta *0 eta) *1(mu *0 eta) *1 mu}}
	\ar@3 @/^9ex/ [rrrr] ^-{\twocell{rho}} _-{}="s1"
	\ar@3 [rr] |-*+[o]{\twocell{alpha}} ^-{}="t1" _-{}="s2"
	\ar@3 @/_9ex/ [rrrr] _-{\twocell{rho}} ^-{}="t2"
&& {\twocell{(eta *0 eta) *1 (eta *0 mu) *1 mu}}
	\ar@3 @/^3ex/ [rr] ^(0.4){\twocell{lambda}} _-{}="s3"
	\ar@3 @/_3ex/ [rr] _(0.4){\twocell{rho}} ^-{}="t3"
&& {\twocell{(eta *0 eta) *1 mu}}
	\ar@3 [rr] ^-{\twocell{rho}}
&& {\twocell{eta}}
\ar@{} "s1";"t1" |-{\twocell{(eta *0 eta) *1 beth}} 
\ar@{} "s2";"t2" |-{\twocell{eta *0 eta} \star_1 \omega_2}
\ar@{} "s3";"t3" |-{\delta}
& = 
& {\twocell{eta *1 rho *1 rho}}
}
}
\]
As a consequence, we have $\cl{s(\delta)}=\cl{t(\delta)}$ in $\Mon\Cat$. Hence, we also have equality
\[
\cl{s(\delta)} \star_2 \cl{\twocell{eta *1 lambda}} \:=\: 
\cl{t(\delta)} \star_2 \cl{\twocell{eta *1 lambda}}.
\]
in $\Mon\Cat$, which relates the source and target of the following diagram:
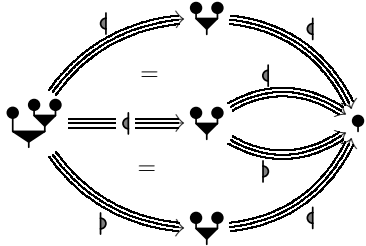
\[
\scalebox{0.8}{
\xymatrix{
&& {\twocell{(eta *0 eta) *1 mu}} 
	\ar@3 @/^3ex/ [drr] ^-{\twocell{lambda}}
\\
{\twocell{(eta *0 ((eta *0 eta) *1 mu)) *1 mu}}
	\ar@3 @/^3ex/ [urr] ^-{\twocell{lambda}} _-{}="s1"
	\ar@3 [rr] |-*+[o]{\twocell{lambda}} 
	\ar@3 @/_3ex/ [drr] _-{\twocell{rho}} ^{}="t2"
&& {\twocell{(eta *0 eta) *1 mu}}
	\ar@3 @/^3ex/ [rr] ^(0.4){\twocell{lambda}}
	\ar@3 @/_3ex/ [rr] _(0.4){\twocell{rho}}
&& {\twocell{eta}}
\\
&& {\twocell{(eta *0 eta) *1 mu}}
	\ar@3 @/_3ex/ [urr] _-{\twocell{lambda}}
\ar@{} "s1";"2,3" |-{=}
\ar@{} "2,3";"t2" |-{=}
}
}
\]
This gives $\cl{s(\omega_3)}=\cl{t(\omega_3)}$ in $\Mon\Cat$, thus concluding the proof.
\end{proof}

\noindent
We can deduce, from this result and Proposition~\ref{generatingSetIar}, that the following two elements form a generating set for the natural system of identities among relations $\Pi(\poly{Mon}_3)$ on the $2$-pro~$\Mon=\cl{\poly{Mon}}_3$ of monoids:
\[
\iar{\tilde{\twocell{aleph}}} \:=\:
\iar{
\twocell{(alpha *0 1) *1 mu}
\star_2
\twocell{(1 *0 mu *0 1) *1 alpha}
\star_2
\twocell{(1 *0 alpha) *1 mu}
\star_2
\left(\twocell{(2 *0 mu)*1 alpha}\right)^{-}
\star_2
\left(\twocell{(mu *0 2)*1 alpha}\right)^{-}
}
\]
\[
\iar{\tilde{\twocell{beth}}} \:=\:
\iar{
\twocell{(1 *0 eta *0 1) *1 alpha}
\star_2
\twocell{(1 *0 lambda) *1 mu}
\star_2
\left(\twocell{(rho *0 1) *1 mu}\right)^{-}
}
\]
From Proposition~\ref{Mon_4HomotopyBase}, we have:

\begin{corollary}[Coherence theorem for monoidal categories,~\cite{MacLane63}]
\label{corollaryMonAspherical}
The $3$-pro $\Mon\Cat$ is aspherical.
\end{corollary}

\section{Coherence in symmetric monoidal categories}
\label{Section:CoherenceSym}

\subsection{Presentations of track \pdf{3}-props}
\label{Subsection:Presentations2Props}

We recall from~\cite{Guiraud04} the following characterisation of $2$-props, derived from a similar result for algebraic theories~\cite{Burroni93}. 

\begin{proposition}
\label{ThmSymmetrique}
A $2$-pro $\P$ is a $2$-prop if and only if it contains a $2$-cell $\tau:2\dfl 2$, represented by~\twocell{tau}, such that the following relations hold:
\begin{itemize}
\item The symmetry relation $\tau\star_1\tau = \id_2$,
\begin{equation}
\label{symmetryRelation}
\twocell{tau *1 tau} \:=\: \twocell{2}
\end{equation}
\item The Yang-Baxter relation $(\tau\star_0 1)\star_1(1\star_0\tau)\star_1(\tau\star_0 1) = (1\star_0\tau)\star_1(\tau\star_0 1)\star_1(1\star_0\tau)$,
\begin{equation}
\label{YangBaxterRelation}
\twocell{(tau *0 1) *1 (1 *0 tau) *1 (tau *0 1)} 
	\:=\: \twocell{(1 *0 tau) *1 (tau *0 1) *1 (1 *0 tau)}
\end{equation}
\item For every $2$-cell $f:m\dfl n$ of $\P$, the left and right naturality relations for $f$,
\[
(f \star_0 1) \star_1 \tau_{n,1} \:=\: \tau_{m,1} \star_1 (1\star_0 f)
\quad\text{and}\quad
(1\star_0 f) \star_1 \tau_{1,n} \:=\: \tau_{1,m} \star_1 (f\star_0 1),
\] 
with the inductively defined notations $\tau_{0,1} = \tau_{1,0} = \id_1$, $\tau_{n+1,1} = (n \star_0 \tau) \star_1 (\tau_{n,1}\star_0 1)$ and $\tau_{1,n+1} = (\tau \star_0 n) \star_1 (1\star_0\tau_{1,n})$. Graphically, we represent $f$ by \twocell{dots *1 f *1 dots}, any $\tau_{n,1}$ by~\twocell{(dots *0 1) *1 tau1 *1 (1 *0 dots)} and any~$\tau_{1,n}$ by~\twocell{(1 *0 dots) *1 tau2 *1 (dots *0 1)}, so that the naturality relations for $f$ are
\begin{equation}
\label{naturalityRelation}
\twocell{(dots *0 1) *1 (f *0 1) *1 (dots *0 1) *1 tau1 *1 (1 *0 dots)} 
	\:=\: \twocell{(dots *0 1) *1 tau1 *1 (1 *0 dots) *1 (1 *0 f) *1 (1 *0 dots)}
\qquad\qquad\text{and}\qquad\qquad
\twocell{(1 *0 dots) *1 (1 *0 f) *1 (1 *0 dots) *1 tau2 *1 (dots *0 1)} 
	\:=\: \twocell{(1 *0 dots) *1 tau2 *1 (dots *0 1) *1 (f *0 1) *1 (dots *0 1)}
\end{equation}
\end{itemize}
\end{proposition}

\subsubsection{The \pdf{2}-prop of permutations}

The initial $2$-prop is the $2$-prop of permutations, denoted by $\Perm$, whose $2$-cells from $n$ to $n$ are the permutations of $\ens{1,\dots,n}$ and with no $2$-cell from $m$ to $n$ if $m\neq n$. The $2$-prop $\Perm$ is presented by the $3$-polygraph with one $2$-cell $\twocell{tau}$ and two $3$-cells, corresponding to the symmetry relation~\eqref{symmetryRelation} and the Yang-Baxter relation~\eqref{YangBaxterRelation}:
\[
\twocell{tau *1 tau} 
	\:\tfl\: \twocell{2} 
\qquad\text{and}\qquad
\twocell{(tau *0 1) *1 (1 *0 tau) *1 (tau *0 1)} 
	\:\tfl\: \twocell{(1 *0 tau) *1 (tau *0 1) *1 (1 *0 tau)}
\]
There exists an isomorphism between the category of small categories and functors and the category $\Alg(\Perm)$. The correspondence between a category $\Cr$ and a $\Perm$-algebra $\Ar : \Perm \fl \Cat$ is given by
\[
\Ar(1)\:=\: \Cr 
\qquad\qquad\text{and}\qquad\qquad
\Ar(\twocell{tau}) \:=\: T_{\Cr,\Cr},
\] 
where $T_{\Cr,\Cr}$ is the endofunctor of $\Cr\times\Cr$ sending $(x,y)$ to $(y,x)$. 

\subsubsection{Presentations of \pdf{2}-props}

Let $\Sigma$ be a $2$-polygraph with one $0$-cell and one $1$-cell. We denote by~$S\Sigma$ the $3$-polygraph obtained from $\Sigma$ by adjoining a $2$-cell $\twocell{tau}:2\dfl 2$ and the following $3$-cells:
\begin{itemize}
\item The symmetry $3$-cell and the Yang-Baxter $3$-cell, as in the $2$-prop $\Perm$.
\item Two $3$-cells for every $2$-cell $f=\twocell{dots *1 f *1 dots}$ of $\Sigma$, corresponding to the naturality relations for $f$:
\[
\twocell{(dots *0 1) *1 (f *0 1) *1 (dots *0 1) *1 tau1 *1 (1 *0 dots)} 
	\:\tfl\: \twocell{(dots *0 1) *1 tau1 *1 (1 *0 dots) *1 (1 *0 f) *1 (1 *0 dots)}
\qquad\qquad\text{and}\qquad\qquad
\twocell{(1 *0 dots) *1 (1 *0 f) *1 (1 *0 dots) *1 tau2 *1 (dots *0 1)} 
	\:\tfl\: \twocell{(1 *0 dots) *1 tau2 *1 (dots *0 1) *1 (f *0 1) *1 (dots *0 1)}
\]
\end{itemize}
The \emph{free $2$-prop generated by $\Sigma$} is the $2$-category, denoted by $\Sigma^S$, presented by the $3$-polygraph $S\Sigma$.  

Let $\P$ be a $2$-prop. A \emph{presentation of $\P$} is a pair $(\Sigma_2,\Sigma_3)$, made of a $2$-polygraph $\Sigma_2$ with one $0$-cell and one $1$-cell and a cellular extension $\Sigma_3$ of the free $2$-prop $\Sigma_2^S$, such that  
\[
\P \:\simeq\: \Sigma_2^S/\Sigma_3.
\]

\begin{proposition}
A $3$-pro $\P$ is a $3$-prop if and only if it contains a $2$-cell $\tau:2\dfl 2$ such that the following relations hold:
\begin{itemize}
\item The symmetry relation~\eqref{symmetryRelation} and the Yang-Baxter relation~\eqref{YangBaxterRelation}.
\item The naturality relations~\eqref{naturalityRelation} for every $2$-cell of $\P$.
\item For every $3$-cell $A:f\tfl g:m\dfl n$, the left and right naturality relations for $A$:
\[
(A \star_0 1) \star_1 \tau_{n,1} \:=\: \tau_{m,1} \star_1 (1\star_0 A)
\qquad\text{and}\qquad
(1\star_0 A) \star_1 \tau_{1,n} \:=\: \tau_{1,m} \star_1 (A\star_0 1).
\]
Graphically, we represent $f$ by \twocell{dots *1 f *1 dots} and $g$ by \twocell{dots *1 g *1 dots}, so that the naturality relations for $A$ are
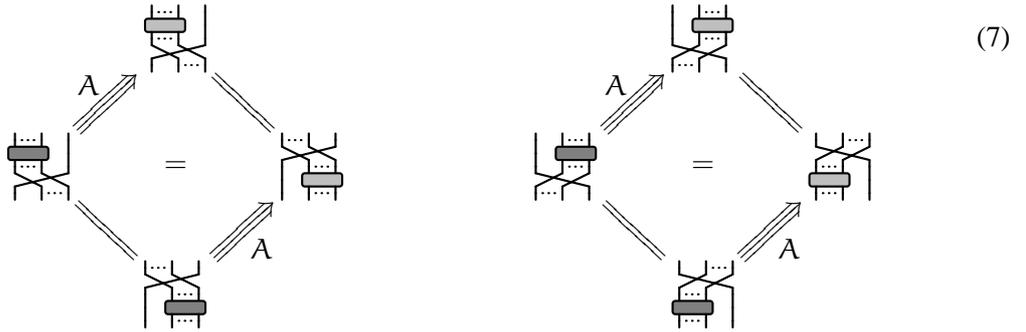
\begin{equation}
\label{naturalityRelation3cells}
\xymatrix@R=1.5em@C=2em{
& {\twocell{(dots *0 1) *1 (g *0 1) *1 (dots *0 1) *1 tau1 *1 (1 *0 dots)}}
	\ar@{=} [dr]
\\
{\twocell{(dots *0 1) *1 (f *0 1) *1 (dots *0 1) *1 tau1 *1 (1 *0 dots)}}
	\ar@3 [ur] ^-{A}
	\ar@{=} [dr]
	\ar@{} [rr] |-{=}
&& {\twocell{(dots *0 1) *1 tau1 *1 (1 *0 dots) *1 (1 *0 g) *1 (1 *0 dots)}}
\\
& {\twocell{(dots *0 1) *1 tau1 *1 (1 *0 dots) *1 (1 *0 f) *1 (1 *0 dots)}}
	\ar@3 [ur] _-{A }
}
\qquad\qquad\qquad
\xymatrix@R=1.5em@C=2em{
& {\twocell{(1 *0 dots) *1 (1 *0 g) *1 (1 *0 dots) *1 tau2 *1 (dots *0 1)}}
	\ar@{=} [dr]
\\
{\twocell{(1 *0 dots) *1 (1 *0 f) *1 (1 *0 dots) *1 tau2 *1 (dots *0 1)}}
	\ar@3 [ur] ^-{A}
	\ar@{=} [dr]
	\ar@{} [rr] |-{=}
&& {\twocell{(1 *0 dots) *1 tau2 *1 (dots *0 1) *1 (g *0 1) *1 (dots *0 1)}}
\\
& {\twocell{(1 *0 dots) *1 tau2 *1 (dots *0 1) *1 (f *0 1) *1 (dots *0 1)}}
	\ar@3 [ur] _-{A }
}
\end{equation}
\end{itemize}
\end{proposition}

\begin{proof}
This is an immediate extension of Proposition~\ref{ThmSymmetrique}.   
\end{proof}

\subsubsection{Presentations of track \pdf{3}-props}

Let $\Sigma$ be a presentation of a $2$-prop. We denote by~$S\Sigma$ the $4$-polygraph obtained from the $3$-polygraph $S\Sigma_2$ by adjoining the $3$-cells of $\Sigma_3$ and a cellular extension~$\Sigma_4$ made of the following two $4$-cells for each $3$-cell $A$ of $\Sigma_3$, corresponding to the naturality relations~\eqref{naturalityRelation3cells} for $A$:
\[
\xymatrix@R=1.5em@C=2em{
& {\twocell{(dots *0 1) *1 (g *0 1) *1 (dots *0 1) *1 tau1 *1 (1 *0 dots)}}
	\ar@3 [dr]
	\ar@4 []!<0pt,-30pt>;[dd]!<0pt,30pt> 
\\
{\twocell{(dots *0 1) *1 (f *0 1) *1 (dots *0 1) *1 tau1 *1 (1 *0 dots)}}
	\ar@3 [ur] ^-{A}
	\ar@3 [dr]
&& {\twocell{(dots *0 1) *1 tau1 *1 (1 *0 dots) *1 (1 *0 g) *1 (1 *0 dots)}}
\\
& {\twocell{(dots *0 1) *1 tau1 *1 (1 *0 dots) *1 (1 *0 f) *1 (1 *0 dots)}}
	\ar@3 [ur] _-{A}
}
\qquad\qquad\qquad
\xymatrix@R=1.5em@C=2em{
& {\twocell{(1 *0 dots) *1 (1 *0 g) *1 (1 *0 dots) *1 tau2 *1 (dots *0 1)}}
	\ar@3 [dr]
	\ar@4 []!<0pt,-30pt>;[dd]!<0pt,30pt> 
\\
{\twocell{(1 *0 dots) *1 (1 *0 f) *1 (1 *0 dots) *1 tau2 *1 (dots *0 1)}}
	\ar@3 [ur] ^-{A}
	\ar@3 [dr]
&& {\twocell{(1 *0 dots) *1 tau2 *1 (dots *0 1) *1 (g *0 1) *1 (dots *0 1)}}
\\
& {\twocell{(1 *0 dots) *1 tau2 *1 (dots *0 1) *1 (f *0 1) *1 (dots *0 1)}}
	\ar@3 [ur] _-{A}
}
\]
The \emph{free track $3$-prop generated by $\Sigma$} is the track $3$-category, denoted by $\Sigma^S$, given by: 
\[
\Sigma^S \:=\: \Sigma_2^S(\Sigma_3) /\Sigma_4.
\]
Let $\P$ be a track $3$-prop. A \emph{presentation of $\P$} is a pair $(\Sigma_3,\Sigma_4)$, where $\Sigma_3$ is a presentation of a $2$-prop and~$\Sigma_4$ is a cellular extension of the free track $3$-prop $\Sigma_3^S$, such that  
\[
\P \:\simeq\: \Sigma_3^S/\Sigma_4.
\]
To summarize, a presentation of $\P$ yields a diagram which is similar to the one corresponding to the inductive construction of a $4$-polygraph, see Section~\ref{Polygraphs}:
\[
\xymatrix{
\ens{0} 
     \ar@{=}[d] 
&& \Nb 
&& \Sigma_2^S
&& \Sigma_3^S
\\
\ens{0}
&& \ens{1}
     \ar@<0.5ex>[llu]
     \ar@<-0.5ex>[llu]
     \ar@{>->}[u]
&& \Sigma_2
     \ar@<0.5ex>[llu]
     \ar@<-0.5ex>[llu]
     \ar@{>->}[u]
&& \Sigma_3
     \ar@<0.5ex>[llu]
     \ar@<-0.5ex>[llu]
     \ar@{>->}[u]
&& \Sigma_4
     \ar@<0.5ex>[llu]
     \ar@<-0.5ex>[llu]
}
\]

\subsection{Convergent presentations of algebraic track \pdf{3}-props and asphericity}

\subsubsection{Convergent presentations of algebraic track \pdf{3}-props}
\label{algebraicPresentationProp}

A presentation $\Sigma$ of a track $3$-prop is \emph{convergent} when the $3$-polygraph $S\Sigma$ is convergent. A presentation $\Sigma$ of a $2$-prop (resp. track $3$-prop) is \emph{algebraic} when every $2$-cell (resp. every $2$-cell and every $3$-cell) of $\Sigma$ has $1$-target equal to the generating $1$-cell $1$. A track $3$-prop is \emph{algebraic} when it admits an algebraic presentation. 

\subsubsection{Classification of critical branchings}

Let $\Sigma$ be an algebraic presentation of a $2$-prop $\P$. We recall from~\cite{Guiraud04,Guiraud06jpaa} that the critical branchings of the $3$-polygraph $S\Sigma$ are classified as follows:
\begin{enumerate}
\item Five critical branchings generated by the symmetry and Yang-Baxter $3$-cells, whose sources are:
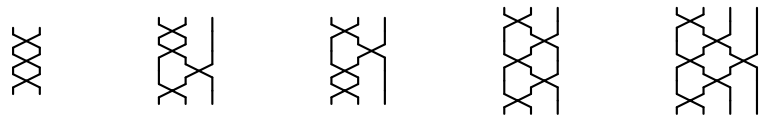
\[
\twocell{tau *1 tau *1 tau}
\qquad\qquad \twocell{(tau *0 1) *1 (tau *0 1) *1 (1 *0 tau) *1 (tau *0 1)}
\qquad\qquad \twocell{(tau *0 1) *1 (1 *0 tau) *1 (tau *0 1) *1 (tau *0 1)}
\qquad\qquad \twocell{(tau *0 1) *1 (1 *0 tau) *1 (tau *0 1) *1 (1 *0 tau) *1 (tau *0 1)}
\qquad\qquad \twocell{(tau *0 2) *1 (1 *0 tau *0 1) *1 (tau *0 tau) *1 (1 *0 tau *0 1) *1 (tau *0 2)}
\]
\item For every $2$-cell $\phi=\twocell{dots *1 phi}$ of $\Sigma$, five critical branchings, generated, on the one hand, by the naturality $3$-cells for $\phi$ and, on the other hand, by the symmetry and Yang-Baxter $3$-cells:
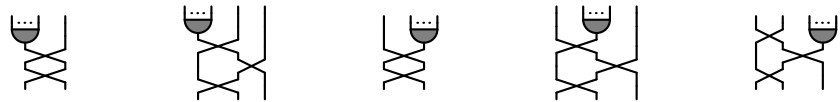
\[
\twocell{(dots *0 1) *1 (phi *0 1) *1 tau *1 tau}
\qquad\qquad \twocell{(dots *0 2) *1 (phi *0 2) *1 (tau *0 1) *1 (1 *0 tau) *1 (tau *0 1)}
\qquad\qquad \twocell{(1 *0 dots) *1 (1 *0 phi) *1 tau *1 tau}
\qquad\qquad \twocell{(1 *0 dots *0 1) *1 (1 *0 phi *0 1) *1 (tau *0 1) *1 (1 *0 tau) *1 (tau *0 1)}
\qquad\qquad \twocell{(2 *0 dots) *1 (tau *0 phi) *1 (1 *0 tau) *1 (tau *0 1)}
\]
\item For every pair $(\phi,\psi)$ of $2$-cells of $\Sigma$, one critical branching generated by the left naturality $3$-cell of $\phi=\twocell{dots *1 phi}$ and the right naturality $3$-cell of $\psi=\twocell{dots *1 psi}$:
\[
\twocell{(dots *0 dots) *1 (phi *0 psi) *1 tau}
\]
\item For every algebraic $3$-cell $\alpha:f\tfl g$ of $\Sigma$, two critical branchings generated by $\alpha$ and the naturality $3$-cells for $f=\twocell{dots *1 phi}$:
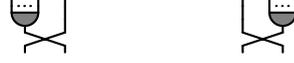
\[
\twocell{(dots *0 1) *1 (phi *0 1) *1 tau}
\qquad\qquad\qquad
\twocell{(1 *0 dots) *1 (1 *0 phi) *1 tau}
\] 
\item The other critical branchings, called the \emph{proper critical branchings of $\Sigma$}.
\end{enumerate}
All of the critical branchings of the first three families are confluent and their confluence diagrams are sent to commutative diagrams by the canonical projection $\pi:\tck{S\Sigma}\fl \Sigma^S$. The critical branchings of the fourth family are confluent and their confluence  diagrams are sent to $3$-spheres which are the boundaries of naturality $4$-cells. 

A \emph{basis of proper confluences of $\Sigma$} is a cellular extension of the free track $3$-category $\tck{\Sigma}$ that contains, for each proper critical branching $b$ of $\Sigma$, one $4$-cell $\omega_b:A\qfl B$, where the $3$-sphere $(A,B)$ is a confluence diagram for $b$. We assume that, when $\Sigma$ is a convergent $3$-polygraph, we have chosen a basis of proper confluences, which we denote by $\Gamma_{\Sigma}$. 

\begin{lemma}
Let $\Sigma$ be an algebraic convergent presentation of a $2$-prop $\P$. Then the image $\pi(\Gamma_{\Sigma})$ of the cellular extension $\Gamma_{\Sigma}$ through the canonical projection $\pi:\tck{S\Sigma}\fl\Sigma^S$ is a homotopy basis of $\P$. 
\end{lemma}

\begin{theorem}
\label{TheoremPropAspherical}
If a track $3$-prop $\P$ admits an algebraic convergent presentation $(\Sigma_3,\Sigma_4)$ such that $\Sigma_4$ is Tietze-equivalent to $\pi(\Gamma_{\Sigma_3})$, then $\P$ is aspherical. 
\end{theorem}

\subsection{Application to symmetric monoidal categories}
\label{Subsection:ApplicationCoherenceSym}

A \emph{symmetric monoidal category} is a monoidal category $(\Cr,\tens,e,\alpha,\lambda,\rho)$ equipped with a natural isomorphism 
\[
\tau_{x,y} : x\tens y \longrightarrow y\tens x,
\]
called the \emph{symmetry} and such that the following two diagrams commute in $\Cr$:   
\begin{equation}
\label{diagramSymmetry}
\scalebox{0.9}{
\raisebox{-20pt}{
\xymatrix@C=1em{
& y\tens x
	\ar[dr] ^-{\tau}
\\
x\tens y 
	\ar[ur] ^-{\tau}
	\ar@{=}[rr] ^-{}="1"
&& x\tens y
	\ar@{} "1,2";"1" |-{\copyright}
}
}
\qquad
\xymatrix@C=1em{
& x\tens (y\tens z) 
	\ar[rr] ^-{\tau}
&& (y\tens z)\tens x
	\ar[dr] ^-{\alpha}
\\
(x\tens y) \tens z 
	\ar[ur] ^-{\alpha}
	\ar[dr] _-{\tau}
	\ar@{} [rrrr] |-{\copyright}
&&&& y \tens (z\tens x)
\\
& (y\tens x)\tens z 
	\ar[rr] _-{\alpha}
&& y\tens(x\tens z)
	\ar[ur] _-{\tau}
}
}
\end{equation}
A symmetric monoidal functor from $\Cr$ to $\Dr$ is a monoidal functor $(F,\upphi,\iota)$ such that the following diagram commutes in~$\Dr$:
\begin{equation}
\label{symmMonoidalFunctorCommDiag}
\scalebox{0.9}{
\xymatrix{
{ Fx\tens Fy }
	\ar[r] ^-{\tau}
        \ar[d] _-{\upphi}
        \ar@{} [dr] |-{\copyright}
&
{Fy\tens Fx}
        \ar[d] ^-{\upphi}
\\
{F(x\tens y)}
	\ar[r] _-{F\tau}
&
{F(y\tens x)}
}
}
\end{equation}

\subsubsection{The track \pdf{3}-prop of symmetric monoidal categories}
\label{2PropSymCat}

Let $\Sym\Cat$ be the track $3$-prop presented by $\poly{Sym}$ given as follows:
\begin{itemize}
\item $\poly{Sym}_2$ is the $2$-polygraph $\poly{Mon}_2$, containing two $2$-cells~\twocell{mu} and~\twocell{eta}.
\item $\poly{Sym}_3$ is the cellular extension of the free $2$-prop $\poly{Sym}_2^S$ generated by $\poly{Sym}_2$ containing the three $3$-cells of $\poly{Mon}_3$ 
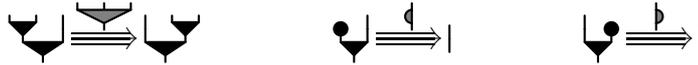
\[
\xymatrix{
{\twocell{ (mu *0 1) *1 mu }}
	\ar@3[r] ^-{\twocell{alpha}}
& {\twocell{ (1 *0 mu) *1 mu } }
}
\qquad\qquad
\xymatrix{
{\twocell{ (eta *0 1) *1 mu }}
	\ar@3[r] ^-{\twocell{lambda}}
& {\twocell{1}}
}  
\qquad\qquad
\xymatrix{
{\twocell{ (1 *0 eta) *1 mu }}
	\ar@3[r] ^-{\twocell{rho}}
& {\twocell{1}}
}
\]
plus the following extra $3$-cell:
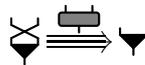
\[
\xymatrix{
{\twocell{tau *1 mu}}
	\ar@3[r] ^-{\twocell{beta}}
& {\twocell{mu}}
}
\]
\item $\poly{Sym}_4$ is the cellular extension of the free $3$-prop $\poly{Sym}_3^S$ generated by $\poly{Sym}_3$ containing the two $4$-cells of $\poly{Mon}_4$ 
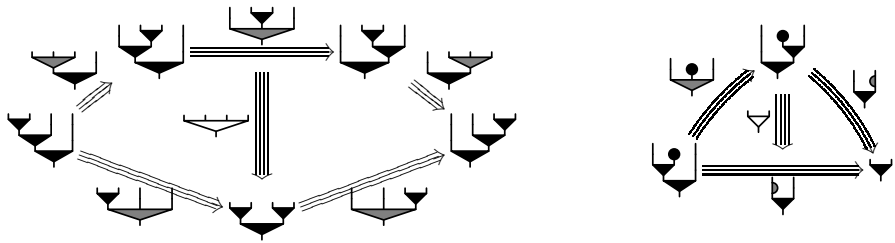
\[
\scalebox{0.8}{
\xymatrix@R=10pt@C=15pt{
&{ \twocell{ (1 *0 mu *0 1) *1 (mu *0 1) *1 mu } }
	\ar@3 [rr] 
		^-{\twocell{(1 *0 mu *0 1) *1 alpha}} 
		_-{}="1" 
&& { \twocell{ (1 *0 mu *0 1) *1 (1 *0 mu) *1 mu } }
	\ar@3 [dr] ^-{\twocell{(1 *0 alpha) *1 mu} }
\\
{ \twocell{ (mu *0 2) *1 (mu *0 1) *1 mu } }
	\ar@3 [ur] ^-{\twocell{(alpha *0 1) *1 mu}}
	\ar@3 [drr] _-{\twocell{(mu *0 2)*1 alpha}}
&&&& { \twocell{ (2 *0 mu) *1 (1 *0 mu) *1 mu } }
\\
&& { \twocell{ (mu *0 mu) *1 mu } }
	\ar@3 [urr] _-{\twocell{(2 *0 mu)*1 alpha}}
	\ar@4 "1"!<0pt,-10pt>;[]!<0pt,20pt> _-*+{\twocell{aleph}}
}
}
\qquad\qquad
\scalebox{0.8}{
\xymatrix{
& { \twocell{ ( 1 *0 eta *0 1) *1 (1 *0 mu) *1 mu } }
	\ar@3@/^/ [dr] ^-{\twocell{(1 *0 lambda) *1 mu}}
\\
{ \twocell{ ( 1 *0 eta *0 1) *1 (mu *0 1) *1 mu } }
	\ar@3@/^/ [ur] ^-{\twocell{(1 *0 eta *0 1) *1 alpha}}
	\ar@3 [rr] _-{\twocell{(rho *0 1) *1 mu}} ^-{}="1"
&& { \twocell{mu} }
\ar@4 "1,2"!<0pt,-20pt>;"1"!<0pt,10pt> _-*+{\twocell{beth}}
}
}
\]
plus the following two extra $4$-cells:
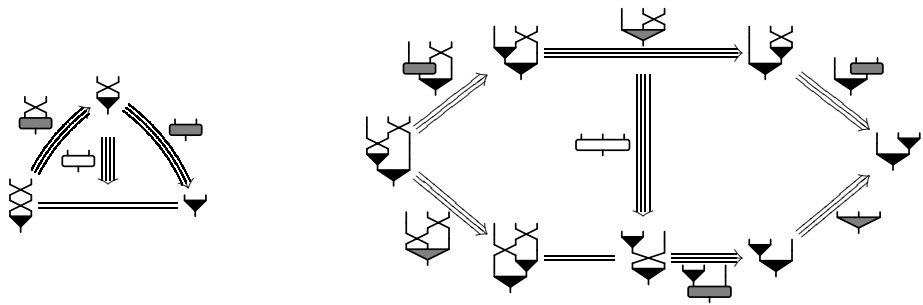
\[
\scalebox{0.8}{
\raisebox{-20pt}{
\xymatrix{
& {\twocell{tau *1 mu}}
	\ar@3@/^/ [dr] ^-{\twocell{beta}}
\\
{\twocell{tau *1 tau *1 mu}}
	\ar@3@/^/ [ur] ^-{\twocell{tau *1 beta}}
	\ar@{=}[rr] ^-{}="1"
&& {\twocell{mu}}
	\ar@4 "1,2"!<0pt,-20pt>;"1"!<0pt,10pt> _-*+{\twocell{gimmel}}
}
}
\qquad\qquad\qquad
\xymatrix@R=1em@C=3em{
& {\twocell{(1 *0 tau) *1 (mu *0 1) *1 mu}}
	\ar@3[rr] ^-{\twocell{(1 *0 tau) *1 alpha}} _-{}="1"
&& {\twocell{(1 *0 tau) *1 (1 *0 mu) *1 mu}}
	\ar@3[dr] ^-{\twocell{(1 *0 beta) *1 mu}}
\\
{\twocell{(1 *0 tau) *1 (tau *0 1) *1 (mu *0 1) *1 mu}}
	\ar@3[ur] ^-{\twocell{(1 *0 tau) *1 (beta *0 1) *1 mu}}
	\ar@3[dr] _-{\twocell{(1 *0 tau) *1 (tau *0 1) *1 alpha}}
&&&& {\twocell{(1 *0 mu) *1 mu}}
\\
& {\twocell{(1 *0 tau) *1 (tau *0 1) *1 (1 *0 mu) *1 mu}}
	\ar@{=} [r] 
& {\twocell{(mu *0 1) *1 tau *1 mu}}
	\ar@3[r] _-{\twocell{(mu *0 1) *1 beta}}
& {\twocell{(mu *0 1) *1 mu}}
	\ar@3[ur] _-{\twocell{alpha}}
\ar@4 "1"!<0pt,-10pt>;"3,3"!<0pt,20pt> _-*+{\twocell{daleth}}
}
}
\]
\end{itemize}

\begin{lemma}
The category of small symmetric monoidal categories and symmetric monoidal functors is isomorphic to the category $\Alg(\Sym\Cat)$.
\end{lemma}

\begin{proof}
Given a symmetric monoidal category $(\Cr,\tens,e,\alpha,\lambda,\rho,\tau)$, the correspondence with a $\Sym\Cat$-algebra~$\Ar$ is given by~\eqref{isoAlgMonEq} for the monoidal underlying structure and by 
\[
\Ar(\twocell{beta}) \:=\: \tau
\]
for the symmetry. The two commutative diagrams of a monoidal category correspond to $\Ar(\twocell{aleph})$ and $\Ar(\twocell{beth})$ and the commutative diagrams~\eqref{diagramSymmetry} correspond to $\Ar(\twocell{gimmel})$ and $\Ar(\twocell{daleth})$. 

The correspondence of a symmetric monoidal functor $(F,\upphi,\iota)$ with a morphism $\Psi$ between the associated $\Sym$-algebras is given by: 
\[
\Psi_{\twocell{1}} \:=\: F,
\qquad
\Psi_{\twocell{mu}} \:=\: \upphi,
\qquad
\Psi_{\twocell{eta}} \:=\: \iota.
\]
The relation~\eqref{symmMonoidalFunctorCommDiag} corresponds to the properties of the morphism $\Psi$. 
\end{proof}

\subsubsection{A convergent presentation of \pdf{\Sym\Cat}}

We define $\poly{Sym}'$ as the presentation $\poly{Sym}$ of $\Sym\Cat$, extended with one $3$-cell
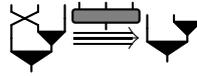
\[
\raisebox{-3em}{
\xymatrix{
{\twocell{(tau *0 1) *1 (1 *0 mu) *1 mu}}
	\ar@3[r] ^-{\twocell{gamma}}
& {\twocell{(1 *0 mu) *1 mu}}
}
}
\]
and the following $4$-cell:
\[
\scalebox{0.8}{
\xymatrix@R=1em@C=2em{
& {\twocell{(mu *0 1) *1 mu}}
	\ar@3 [dr] ^-{\twocell{alpha}}
	\ar@4 []!<0pt,-20pt>;[dd]!<0pt,20pt> _-*+{\omega}
\\
{\twocell{(tau *0 1) *1 (mu *0 1) *1 mu}}
	\ar@3 [ur] ^-{\twocell{(beta *0 1) *1 mu}}
	\ar@3 [dr] _-{\twocell{(tau *0 1) *1 alpha}}
&& {\twocell{(1 *0 mu) *1 mu}}
\\
& {\twocell{(tau *0 1) *1 (1 *0 mu) *1 mu}}
	\ar@3 [ur] _-{\twocell{gamma}}
}
}
\]

\begin{lemma}
The track $3$-prop $\Sym\Cat$ is presented by $\poly{Sym}'$.
\end{lemma}

\begin{proof}
The $4$-cell $\omega$ induces the relation 
\[
\cl{\twocell{gamma}}
	\:=\: 
	\left(\cl{\twocell{(tau *0 1) *1 alpha}}\right)^- 
	\star_2 \cl{\twocell{(beta *0 1) *1 mu}} 
	\star_2 \cl{\twocell{alpha}}
\]
in the quotient track $3$-prop $\left( \poly{Sym}'_3 \right)^S/\poly{Sym}_4'$. As a consequence, it is isomorphic to the quotient track $3$-prop $\Sym\Cat=\left( \poly{Sym}_3 \right)^S/\poly{Sym}_4$.
\end{proof}

\begin{proposition}
The $3$-polygraph $S(\poly{Sym}'_3)$ is convergent and the cellular extension $\poly{Sym}'_4$ is Tietze-equivalent to $\pi(\Gamma_{S(\poly{Sym}'_3)})$.
\end{proposition}

\begin{proof}
The convergence of the $3$-polygraph $S(\poly{Sym}'_3)$ is proved in~\cite{Guiraud06jpaa}. The image through the canonical projection $\pi:\tck{S(\poly{Sym}'_3)}\fl(\poly{Sym}'_3)^S$ of the cellular extension $\Gamma_{S(\poly{Sym}'_3)}$ has ten $4$-cells. Indeed, it contains the images of four $4$-cells of $\poly{Sym}_4$
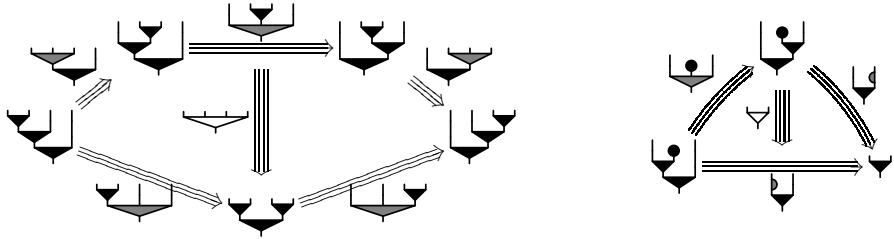
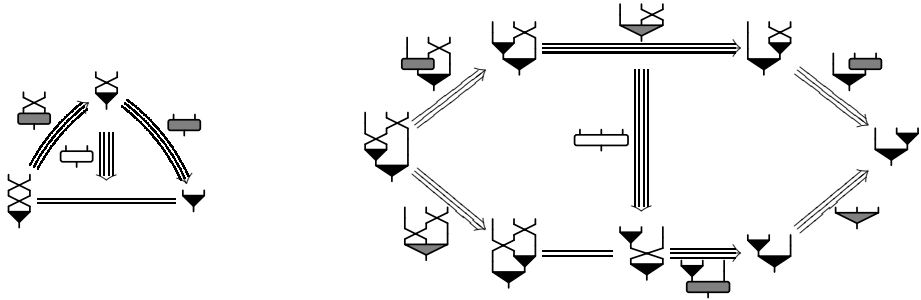
\[
\scalebox{0.8}{
\xymatrix@R=10pt@C=15pt{
&{ \twocell{ (1 *0 mu *0 1) *1 (mu *0 1) *1 mu } }
	\ar@3 [rr] 
		^-{\twocell{(1 *0 mu *0 1) *1 alpha}} 
		_-{}="1" 
&& { \twocell{ (1 *0 mu *0 1) *1 (1 *0 mu) *1 mu } }
	\ar@3 [dr] ^-{\twocell{(1 *0 alpha) *1 mu} }
\\
{ \twocell{ (mu *0 2) *1 (mu *0 1) *1 mu } }
	\ar@3 [ur] ^-{\twocell{(alpha *0 1) *1 mu}}
	\ar@3 [drr] _-{\twocell{(mu *0 2)*1 alpha}}
&&&& { \twocell{ (2 *0 mu) *1 (1 *0 mu) *1 mu } }
\\
&& { \twocell{ (mu *0 mu) *1 mu } }
	\ar@3 [urr] _-{\twocell{(2 *0 mu)*1 alpha}}
	\ar@4 "1"!<0pt,-10pt>;[]!<0pt,20pt> _-*+{\twocell{aleph}}
}
}
\qquad\qquad
\scalebox{0.8}{
\xymatrix{
& { \twocell{ ( 1 *0 eta *0 1) *1 (1 *0 mu) *1 mu } }
	\ar@3@/^/ [dr] ^-{\twocell{(1 *0 lambda) *1 mu}}
\\
{ \twocell{ ( 1 *0 eta *0 1) *1 (mu *0 1) *1 mu } }
	\ar@3@/^/ [ur] ^-{\twocell{(1 *0 eta *0 1) *1 alpha}}
	\ar@3 [rr] _-{\twocell{(rho *0 1) *1 mu}} ^-{}="1"
&& { \twocell{mu} }
\ar@4 "1,2"!<0pt,-20pt>;"1"!<0pt,10pt> _-*+{\twocell{beth}}
}
}
\]
and
\[
\scalebox{0.8}{
\raisebox{-20pt}{
\xymatrix{
& {\twocell{tau *1 mu}}
	\ar@3@/^/ [dr] ^-{\twocell{beta}}
\\
{\twocell{tau *1 tau *1 mu}}
	\ar@3@/^/ [ur] ^-{\twocell{tau *1 beta}}
	\ar@{=}[rr] ^-{}="1"
&& {\twocell{mu}}
	\ar@4 "1,2"!<0pt,-20pt>;"1"!<0pt,10pt> _-*+{\twocell{gimmel}}
}
}
\qquad\qquad\qquad
\xymatrix@R=1em@C=3em{
& {\twocell{(1 *0 tau) *1 (mu *0 1) *1 mu}}
	\ar@3[rr] ^-{\twocell{(1 *0 tau) *1 alpha}} _-{}="1"
&& {\twocell{(1 *0 tau) *1 (1 *0 mu) *1 mu}}
	\ar@3[dr] ^-{\twocell{(1 *0 beta) *1 mu}}
\\
{\twocell{(1 *0 tau) *1 (tau *0 1) *1 (mu *0 1) *1 mu}}
	\ar@3[ur] ^-{\twocell{(1 *0 tau) *1 (beta *0 1) *1 mu}}
	\ar@3[dr] _-{\twocell{(1 *0 tau) *1 (tau *0 1) *1 alpha}}
&&&& {\twocell{(1 *0 mu) *1 mu}}
\\
& {\twocell{(1 *0 tau) *1 (tau *0 1) *1 (1 *0 mu) *1 mu}}
	\ar@{=} [r] 
& {\twocell{(mu *0 1) *1 tau *1 mu}}
	\ar@3[r] _-{\twocell{(mu *0 1) *1 beta}}
& {\twocell{(mu *0 1) *1 mu}}
	\ar@3[ur] _-{\twocell{alpha}}
\ar@4 "1"!<0pt,-10pt>;"3,3"!<0pt,20pt> _-*+{\twocell{daleth}}
}
}
\]
plus the extra $4$-cell $\omega$ of $\poly{Sym}'_4$
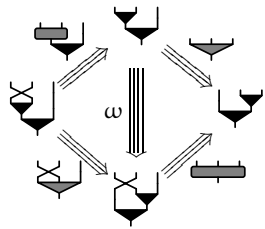
\[
\scalebox{0.8}{
\xymatrix@R=1em@C=2em{
& {\twocell{(mu *0 1) *1 mu}}
	\ar@3 [dr] ^-{\twocell{alpha}}
	\ar@4 []!<0pt,-20pt>;[dd]!<0pt,20pt> _-*+{\omega}
\\
{\twocell{(tau *0 1) *1 (mu *0 1) *1 mu}}
	\ar@3 [ur] ^-{\twocell{(beta *0 1) *1 mu}}
	\ar@3 [dr] _-{\twocell{(tau *0 1) *1 alpha}}
&& {\twocell{(1 *0 mu) *1 mu}}
\\
& {\twocell{(tau *0 1) *1 (1 *0 mu) *1 mu}}
	\ar@3 [ur] _-{\twocell{gamma}}
}
}
\]
and, finally, the following five $4$-cells:
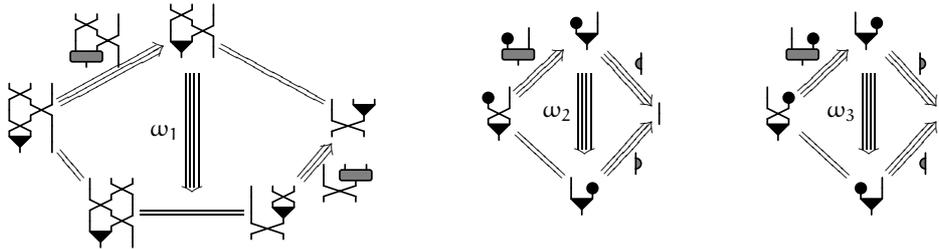
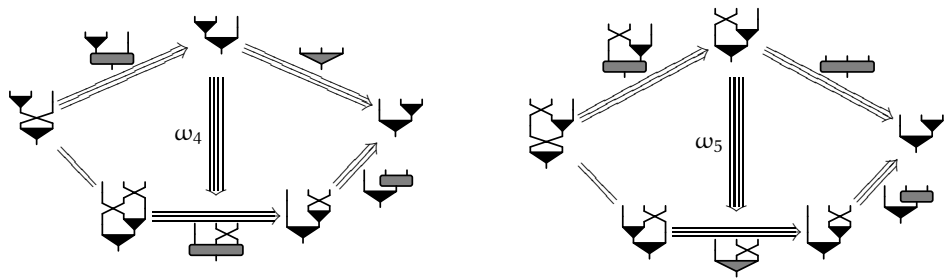
\[
\scalebox{0.8}{
\xymatrix@R=0.5em@C=1em{
&& {\twocell{(tau *0 1) *1 (1 *0 tau) *1 (mu *0 1)}}
	\ar@{=} [drr]
\\
{ \twocell{(tau *0 1) *1 (1 *0 tau) *1 (tau *0 1) *1 (mu *0 1)} }
	\ar@3 [urr] ^-{\twocell{(tau *0 1) *1 (1 *0 tau) *1 (beta *0 1)}}
	\ar@{=} [dr] 
&&&& { \twocell{(1 *0 mu) *1 tau} }
\\
&{ \twocell{(1 *0 tau) *1 (tau *0 1) *1 (1 *0 tau) *1 (mu *0 1)} }
	\ar@{=} [rr] ^-{}="1" 
&& { \twocell{ (1 *0 tau) *1 (1 *0 mu) *1 tau} }
	\ar@3 [ur] _-{\twocell{(1 *0 beta) *1 tau} }
\ar@4 "1,3"!<0pt,-20pt>;"1"!<0pt,10pt> _-*+{\omega_1}
}
\qquad\qquad
\xymatrix@R=1em@C=2em{
& {\twocell{(eta *0 1) *1 mu}}
	\ar@3 [dr] ^-{\twocell{lambda}}
	\ar@4 []!<0pt,-20pt>;[dd]!<0pt,20pt> _-*+{\omega_2}
\\
{\twocell{(eta *0 1) *1 tau *1 mu}}
	\ar@3 [ur] ^-{\twocell{(eta *0 1) *1 beta}}
	\ar@{=} [dr] 
&& {\twocell{1}}
\\
& {\twocell{(1 *0 eta) *1 mu}}
	\ar@3 [ur] _-{\twocell{rho}}
}
\qquad\qquad
\xymatrix@R=1em@C=2em{
& {\twocell{(1 *0 eta) *1 mu}}
	\ar@3 [dr] ^-{\twocell{rho}}
	\ar@4 []!<0pt,-20pt>;[dd]!<0pt,20pt> _-*+{\omega_3}
\\
{\twocell{(1 *0 eta) *1 tau *1 mu}}
	\ar@3 [ur] ^-{\twocell{(1 *0 eta) *1 beta}}
	\ar@{=} [dr] 
&& {\twocell{1}}
\\
& {\twocell{(eta *0 1) *1 mu}}
	\ar@3 [ur] _-{\twocell{lambda}}
}
}
\]
and
\[
\scalebox{0.8}{
\xymatrix@R=1em@C=1.5em{
&& { \twocell{ (mu *0 1) *1 mu } }
	\ar@3 [drr] ^-{ \twocell{alpha} }
\\
{ \twocell{(mu *0 1) *1 tau *1 mu} }
	\ar@3 [urr] ^-{\twocell{ (mu *0 1) *1 beta }}
	\ar@{=} [dr] 
&&&& { \twocell{ (1 *0 mu) *1 mu } }
\\
&{ \twocell{ (1 *0 tau) *1 (tau *0 1) *1 (1 *0 mu) *1 mu } }
	\ar@3 [rr] _-{\twocell{(1 *0 tau) *1 gamma}} ^-{}="1" 
&& { \twocell{ (1 *0 tau) *1 (1 *0 mu) *1 mu } }
	\ar@3 [ur] _-{\twocell{ (1 *0 beta) *1 mu } }
\ar@4 "1,3"!<0pt,-20pt>;"1"!<0pt,10pt> _-*+{ \omega_4 }
}
\qquad\qquad
\xymatrix@R=1em@C=1.5em{
&& { \twocell{ (tau *0 1) *1 (1 *0 mu) *1 mu } }
	\ar@3 [drr] ^-{\twocell{gamma}}
\\
{ \twocell{(tau *0 1) *1 (1 *0 mu) *1 tau *1 mu} }
	\ar@3 [urr] ^-{\twocell{ (tau *0 1) *1 (1 *0 mu) *1 beta }}
	\ar@{=} [dr] 
&&&& { \twocell{ (1 *0 mu) *1 mu } }
\\
&{ \twocell{ (1 *0 tau) *1 (mu *0 1) *1 mu } }
	\ar@3 [rr] 
		_-{\twocell{ (1 *0 tau) *1 alpha }} 
		^-{}="1" 
&& { \twocell{ (1 *0 tau) *1 (1 *0 mu) *1 mu } }
	\ar@3 [ur] _-{\twocell{ (1 *0 beta) *1 mu } }
\ar@4 "1,3"!<0pt,-20pt>;"1"!<0pt,10pt> _-*+{ \omega_5 }
}
}
\]
In order to show that $\poly{Sym}'_4$ is Tietze-equivalent to $\pi(\Gamma_{S(\poly{Sym}'_3)})$, we check that, for each one of the five $4$-cells $\omega_i$, we have the relation $\cl{s(\omega_i)}= \cl{t(\omega_i)}$ in the quotient track $3$-prop $\Sym\Cat$. The projection sends $\omega_1$ to one of the naturality relations for \twocell{beta}. For each one of the other $4$-cells $\omega_i$, with $2\leq i\leq 5$, we consider a $4$-cell $W_i$ of the track $4$-prop $\tck{\poly{Sym}}_3(\poly{Sym}_4)$, built as an instance of the $4$-cell~\twocell{daleth} composed with $2$-cells:
\[
W_2 \:=\: \twocell{(eta *0 eta *0 1) *1 daleth} \;,
\qquad W_3 \:=\: \twocell{(1 *0 eta *0 eta) *1 daleth} \;, 
\qquad W_4 \:=\: \twocell{(mu *0 eta *0 1) *1 daleth} \;,
\qquad W_5 \:=\: \twocell{(tau *0 1) *1 (1 *0 eta *0 mu) *1 daleth} \;.
\]
On the one hand, by definition, the boundary of $W_i$ satisfies the relation $\cl{s(W_i)}=\cl{t(W_i)}$ in the quotient track $3$-prop $\Sym\Cat$. On the other hand, we progressively fill the boundary of $W_i$, as in the case of the track $3$-prop $\Mon\Cat$, with $4$-cells of $\poly{Sym}_4$, plus exchange and naturality relations, until reaching the boundary of the $4$-cell $\omega_i$ (or of $\omega_i^-$), thus yielding the result. 
\end{proof}

\begin{corollary}[Coherence theorem for symmetric monoidal     categories,~\cite{MacLane63}] 
\label{CoherenceTheoremSymmetricMonoidalCategories}
The track $3$-prop $\Sym\Cat$ is aspherical.
\end{corollary}

\section{Coherence for braided monoidal categories}

\subsection{Generalised coherence problem}

A \emph{braided monoidal category} is a monoidal category  $(\Cr, \tens,e,\alpha,\lambda,\rho)$ equipped with a natural isomorphism
\[
\beta_{x,y} : x\tens y \longrightarrow y\tens x,
\]
called the \emph{braiding} and such that the following
diagrams commute in $\Cr$:
\[
\scalebox{0.9}{
\xymatrix@C=1em{
& x\tens (y\tens z) 
	\ar[rr] ^-{\beta}
&& (y\tens z)\tens x
	\ar[dr] ^-{\alpha}
\\
(x\tens y) \tens z 
	\ar[ur] ^-{\alpha}
	\ar[dr] _-{\beta}
	\ar@{} [rrrr] |-{\copyright}
&&&& y \tens (z\tens x)
\\
& (y\tens x)\tens z 
	\ar[rr] _-{\alpha}
&& y\tens(x\tens z)
	\ar[ur] _-{\beta}
}
}
\]  
and
\[
\scalebox{0.9}{
\xymatrix@C=1em{
& x\tens (y\tens z) 
	\ar[rr] ^-{\beta^-}
&& (y\tens z)\tens x
	\ar[dr] ^-{\alpha}
\\
(x\tens y) \tens z 
	\ar[ur] ^-{\alpha}
	\ar[dr] _-{\beta^-}
	\ar@{} [rrrr] |-{\copyright}
&&&& y \tens (z\tens x)
\\
& (y\tens x)\tens z 
	\ar[rr] _-{\alpha}
&& y\tens(x\tens z)
	\ar[ur] _-{\beta^-}
}
}
\]
A braided monoidal functor from $\Cr$ to $\Dr$ is a monoidal functor $(F,\upphi,\iota)$ such that the following diagram commutes in~$\Dr$:
\[
\scalebox{0.9}{
\xymatrix{
{ Fx\tens Fy }
	\ar[r] ^-{\beta}
        \ar[d] _-{\upphi}
        \ar@{} [dr] |-{\copyright}
&
{Fy\tens Fx}
        \ar[d] ^-{\upphi}
\\
{F(x\tens y)}
	\ar[r] _-{F\beta}
&
{F(y\tens x)}
}
}
\]

\subsubsection{Generalised coherence theorem}

Contrary to the case of monoidal and symmetric monoidal categories, we do not have that every diagram commutes in a braided monoidal category. For example, the morphisms $\beta_{x,y}$  and $\beta^{-}_{y,x}$, from $x\tens y$ to $y\tens x$, have no reason to be equal. In fact, they are equal if and only if $\beta$ is a symmetry, hence if and only if all diagrams commute.  

As a consequence, the coherence problem for braided monoidal categories requires a generalised version of the coherence problem we have considered so far. 

\medskip
\noindent
\textsc{The generalised coherence problem:} 

\medskip
\begin{minipage}{14cm}
\emph{Given a track $3$-prop $\P$, decide, for any $3$-sphere $\gamma$ of $\P$, whether or not the diagram $\Ar(\gamma)$ commutes in every $\P$-algebra $\Ar$.}
\end{minipage}

\medskip\smallskip
\noindent
Hence, a solution for the generalised coherence problem is a decision procedure for the equality of $3$-cells of $\P$. For the coherence problems considered so far, this decision procedure  answers yes for every $3$-sphere.  We consider methods to study the generalised coherence theorem of $3$-props and we illustrate those methods on the track $3$-prop of braided monoidal categories.

\subsubsection{The track \pdf{3}-prop of braided monoidal categories}
\label{2PropBrCat}

Let $\Brd\Cat$ be the track $3$-prop with the presentation~$\poly{Br}$ defined
as follows: 
\begin{itemize}
\item The $2$-polygraph $\poly{Br}_2$ is $\poly{Mon}_2$, containing the two $2$-cells \twocell{mu} and \twocell{eta}.
\item The cellular extension $\poly{Br}_3$ of $\poly{Br}_2^S$ has the same four $3$-cells as $\poly{Sym}_3$:
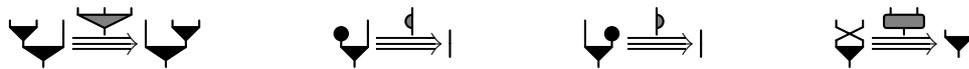
\[
\xymatrix{
{\twocell{ (mu *0 1) *1 mu }}
	\ar@3[r] ^-{\twocell{alpha}}
& {\twocell{ (1 *0 mu) *1 mu } }
}
\qquad\qquad
\xymatrix{
{\twocell{ (eta *0 1) *1 mu }}
	\ar@3[r] ^-{\twocell{lambda}}
& {\twocell{1}}
}  
\qquad\qquad
\xymatrix{
{\twocell{ (1 *0 eta) *1 mu }}
	\ar@3[r] ^-{\twocell{rho}}
& {\twocell{1}}
}
\qquad\qquad
\xymatrix{
{\twocell{tau *1 mu}}
	\ar@3[r] ^-{\twocell{beta}}
& {\twocell{mu}}
}
\]
\item The cellular extension $\poly{Br}_4$ of $\poly{Br}_3^S$ has four $4$-cells, the two $4$-cells of $\poly{Mon}_4$
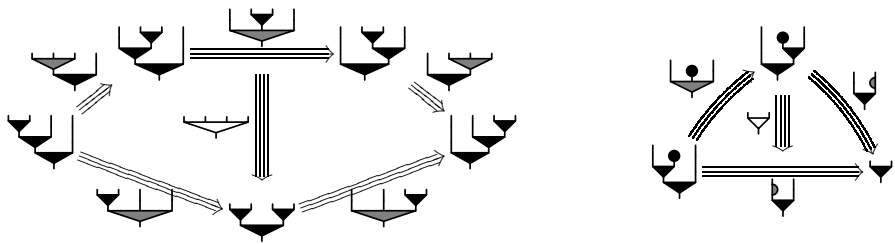
\[
\scalebox{0.8}{
\xymatrix@R=10pt@C=15pt{
&{ \twocell{ (1 *0 mu *0 1) *1 (mu *0 1) *1 mu } }
	\ar@3 [rr] 
		^-{\twocell{(1 *0 mu *0 1) *1 alpha}} 
		_-{}="1" 
&& { \twocell{ (1 *0 mu *0 1) *1 (1 *0 mu) *1 mu } }
	\ar@3 [dr] ^-{\twocell{(1 *0 alpha) *1 mu} }
\\
{ \twocell{ (mu *0 2) *1 (mu *0 1) *1 mu } }
	\ar@3 [ur] ^-{\twocell{(alpha *0 1) *1 mu}}
	\ar@3 [drr] _-{\twocell{(mu *0 2)*1 alpha}}
&&&& { \twocell{ (2 *0 mu) *1 (1 *0 mu) *1 mu } }
\\
&& { \twocell{ (mu *0 mu) *1 mu } }
	\ar@3 [urr] _-{\twocell{(2 *0 mu)*1 alpha}}
	\ar@4 "1"!<0pt,-10pt>;[]!<0pt,20pt> _-*+{\twocell{aleph}}
}
}
\qquad\qquad
\scalebox{0.8}{
\xymatrix{
& { \twocell{ ( 1 *0 eta *0 1) *1 (1 *0 mu) *1 mu } }
	\ar@3@/^/ [dr] ^-{\twocell{(1 *0 lambda) *1 mu}}
\\
{ \twocell{ ( 1 *0 eta *0 1) *1 (mu *0 1) *1 mu } }
	\ar@3@/^/ [ur] ^-{\twocell{(1 *0 eta *0 1) *1 alpha}}
	\ar@3 [rr] _-{\twocell{(rho *0 1) *1 mu}} ^-{}="1"
&& { \twocell{mu} }
\ar@4 "1,2"!<0pt,-20pt>;"1"!<0pt,10pt> _-*+{\twocell{beth}}
}
}
\]
plus the following two $4$-cells:
\[
\scalebox{0.8}{
\xymatrix@R=1em@C=4em{
& {\twocell{(1 *0 tau) *1 (mu *0 1) *1 mu}}
	\ar@3[rr] ^-{\twocell{(1 *0 tau) *1 alpha}} _-{}="1"
&& {\twocell{(1 *0 tau) *1 (1 *0 mu) *1 mu}}
	\ar@3[dr] ^-{\twocell{(1 *0 beta) *1 mu}}
\\
{\twocell{(1 *0 tau) *1 (tau *0 1) *1 (mu *0 1) *1 mu}}
	\ar@3[ur] ^-{\twocell{(1 *0 tau) *1 (beta *0 1) *1 mu}}
	\ar@3[dr] _-{\twocell{(1 *0 tau) *1 (tau *0 1) *1 alpha}}
&&&& {\twocell{(1 *0 mu) *1 mu}}
\\
& {\twocell{(1 *0 tau) *1 (tau *0 1) *1 (1 *0 mu) *1 mu}}
	\ar@{=} [r] 
& {\twocell{(mu *0 1) *1 tau *1 mu}}
	\ar@3[r] _-{\twocell{(mu *0 1) *1 beta}}
& {\twocell{(mu *0 1) *1 mu}}
	\ar@3[ur] _-{\twocell{alpha}}
\ar@4 "1"!<0pt,-10pt>;"3,3"!<0pt,20pt> _-*+{\twocell{daleth}_1}
}
}
\]
and
\[
\scalebox{0.8}{
\xymatrix@R=1em@C=4em{
& {\twocell{(tau *0 1) *1 (mu *0 1) *1 mu}}
	\ar@3[rr] ^-{\twocell{(tau *0 1) *1 alpha}} _-{}="1"
&& {\twocell{(tau *0 1) *1 (1 *0 mu) *1 mu}}
	\ar@3[dr] ^-{{\left(\twocell{(tau *0 1) *1 (1 *0 beta) *1 mu}\right)}^-}
\\
{\twocell{(mu *0 1) *1 mu}}
	\ar@3[ur] ^-{{\left(\twocell{(beta *0 1) *1 mu}\right)}^-}
	\ar@3[dr] _-{\twocell{alpha}}
&&&& {\twocell{(tau *0 1) *1 (1 *0 tau) *1 (1 *0 mu) *1 mu}}
\\
& {\twocell{(1 *0 mu) *1 mu}}
	\ar@3 [r] _{{\left(\twocell{(1 *0 mu) *1 beta}\right)}^-}
& {\twocell{(1 *0 mu) *1 tau *1 mu}}
	\ar@{=} [r]
& {\twocell{(tau *0 1) *1 (1 *0 tau) *1 (mu *0 1) *1 mu}}
	\ar@3[ur] _-{\twocell{(tau *0 1) *1 (1 *0 tau) *1 alpha}}
\ar@4 "1"!<0pt,-10pt>;"3,3"!<0pt,20pt> _-*+{\twocell{daleth}_2}
}
}
\]
\end{itemize}
The correspondence between symmetric monoidal categories and $\Sym\Cat$-algebras can be extended to braided monoidal categories:

\begin{lemma}
The category of small braided monoidal categories and braided monoidal functors is isomorphic to the category $\Alg(\Brd\Cat)$.
\end{lemma}

\subsection{Preservation of coherence by equivalences}

\subsubsection{Equivalence of track \pdf{3}-props}

Let $\P$ and $\Q$ be track $3$-props. A \emph{morphism of track $3$-props} from~$\P$ to $\Q$ is a $3$-functor $F:\P\fl\Q$ which is the identity on $1$-cells, \ie, $F(n)=n$ for every $1$-cell~$n$ in~$\Nb$. If $F,G:\P\fl\Q$ are two morphisms of track $3$-props, a \emph{natural transformation from $F$ to $G$} is a family $\alpha$ of $3$-cells of $\Q$ 
\[
\alpha_f \::\: F(f) \:\tfl\: G(f)
\]
indexed by the $2$-cells of $\P$ and such that, for every $3$-cell $A:f\tfl g$ of $\P$, the following diagram commutes in $\Q$:
\[
\xymatrix{
F(f) 
	\ar@3 [rr] ^-{\alpha_f}
	\ar@3 [d] _-{F(A)}
	\ar@{} [drr] |-{\copyright}
&& G(f)
	\ar@3 [d] ^-{G(A)}
\\
F(g)
	\ar@3 [rr] _-{\alpha_g}
&& G(g).
}
\]
If $F:\P\fl\Q$ is a morphism of track $3$-props, a \emph{quasi-inverse for $F$} is a morphism of track $3$-props $G:\Q\fl\P$ such that there exist natural isomorphisms
\[
GF \:\simeq\: \id_{\P}
\qquad\text{and}\qquad
FG \:\simeq\: \id_{\Q}.
\]
An \emph{equivalence between $\P$ and $\Q$} is a morphism of track $3$-props $F:\P\fl\Q$ that admits a quasi-inverse.

\begin{proposition}
Let $F:\P\fl\Q$ be an equivalence between track $3$-props $\P$ and $\Q$ and let $(A,B)$ be a $3$-sphere of~$\P$. Then $A=B$ if and only if $F(A)=F(B)$. 
\end{proposition}

\begin{proof}
Let $(A,B):f\tfl g$ be a $3$-sphere of $\P$ such that $F(A)=F(B)$. We denote by $G:\Q\fl\P$ a quasi-inverse of $F$ and by $\alpha$ the natural isomorphism from $GF$ to $\id_{\P}$. We have, by definition of $\alpha$, commutative diagrams in $\P$: 
\[
\xymatrix{
GF(f) 
	\ar@3 [rr] ^-{\alpha_f}
	\ar@3 [d] _-{GF(A)}
	\ar@{} [drr] |-{\copyright}
&& f
	\ar@3 [d] ^-{A}
\\
GF(g)
	\ar@3 [rr] _-{\alpha_g}
&& g
}
\qquad\qquad
\xymatrix{
GF(f) 
	\ar@3 [rr] ^-{\alpha_f}
	\ar@3 [d] _-{GF(B)}
	\ar@{} [drr] |-{\copyright}
&& f
	\ar@3 [d] ^-{B}
\\
GF(g)
	\ar@3 [rr] _-{\alpha_g}
&& g
}
\]
By hypothesis, we have $GF(A)=GF(B)$. Thus:
\[
A 
	\:=\: \alpha_f^{-}\star_2 GF(A) \star_2 \alpha_g 
	\:=\: \alpha_f^{-}\star_2 GF(B) \star_2 \alpha_g 
	\:=\: B.
	\qedhere
\]
\end{proof}

\subsection{Preservation of coherence by aspherical quotient}

If $\P$ and $\Q$ are track $3$-props with $\Q\subseteq\P$, we denote by $\P/\Q$ the quotient of $\P$ by the $3$-cells of $\Q$ and by $\pi:\P\fl\P/\Q$ the canonical projection.

\begin{theorem}\label{thmEquivProps}
Let $\P$ and $\Q$ be track $3$-props with $\Q$ aspherical and $\Q\subseteq\P$. Then, for every $3$-sphere $(A,B)$ of $\P$, we have $A=B$ if and only if $\pi(A)=\pi(B)$.
\end{theorem}

\begin{proof}
Let $(f,g)$ be a $2$-sphere of $\P$. Since $\Q$ is aspherical, the $3$-cells of $\P$ from $f$ to $g$ are in bijective correspondence with the $3$-cells of $\P/\Q$ from $\pi(f)$ to $\pi(g)$.
\end{proof}

\noindent
By Corollary \ref{corollaryMonAspherical}, the track $3$-pro(p) $\Mon\Cat$ is aspherical, so that we have:

\begin{corollary}
Let $(A,B)$ be a $3$-sphere of $\Brd\Cat$. Then we have $A=B$ in $\Brd\Cat$ if and only if we have $\pi(A)=\pi(B)$ in $\Brd\Cat/\Mon\Cat$. 
\end{corollary}

\subsection{The initial algebra of an algebraic \pdf{2}-prop}
\label{initialAlgebra}

\subsubsection{Algebraic cells} 

Let $\P$ be an algebraic $2$-prop, with an algebraic presentation $\Sigma$. A $2$-cell $f$ of $\P$ is \emph{purely algebraic} when it is algebraic, \ie, it has target $1$, and it is the image of a $2$-cell of $\Sigma_2^*$ by the canonical projection $\Sigma_2^*\fl\P$, \ie, it contains no $2$-cell \twocell{tau}. 

If $f:n\dfl 1$ is an algebraic $2$-cell of $\P$, then the naturality relations~\eqref{naturalityRelation} satisfied in $\P$ imply that $f$ can be decomposed, in a unique way, as
\[
f \:=\: \sigma_f \star_1 \rep{f},
\]
where $\rep{f}:n\dfl 1$ is a purely algebraic $2$-cell of $\P$ and $\sigma_f:n\dfl n$ is the image of a $2$-cell of $\Perm$ by the canonical inclusion $\Perm\fl\P$, \ie, a $2$-cell of $\P$ written with \twocell{tau} only. If we identify $\sigma_f$ with the corresponding permutation of $\ens{1,\dots,n}$, we have, for every $\P$-algebra $\Ar$, the following relation, for every family $(x_1,\dots,x_n)$ of objects of the category $\Ar(1)$:
\[
\Ar(f) (x_1,\dots,x_n) \:=\: \Ar(\rep{f})(x_{\sigma_f(1)},\dots,x_{\sigma_f(n)}).
\]
The $2$-cell $\rep{f}$ can be identified with the equivalence class of $f$ modulo the congruence generated by  
\[
\sigma \star_1 f \:\approx\: \tau \star_1 f
\]
for any permutations $\sigma$ and $\tau$. Similarly, we denote by $\rep{A}$ the equivalence class of an algebraic $3$-cell $A$ of $\P$ modulo the congruence generated by 
\[
\sigma \star_1 A \:\approx\: \tau \star_1 A
\]
for any permutations $\sigma$ and $\tau$.

\subsubsection{Initial algebras}

Let $\P$ be an algebraic track $3$-prop. The \emph{initial $\P$-algebra} is the $\P$-algebra $\Pr$ defined as follows. The category $\Pr(1)$ is given by: 
\begin{itemize}
\item Its objects are the purely algebraic $2$-cells of $\P$, \ie, the equivalence classes $\rep{f}$, for $f$ any algebraic $2$-cell of $\P$.
\item Its morphisms are the equivalence classes $\rep{A}$ for $A$ any algebraic $3$-cell of $\P$. For such a $3$-cell $A:f\tfl g:n\dfl 1$, the corresponding morphism $\rep{A}$ of $\Pr$ has source $\rep{f}$ and target $\rep{g}$.
\item The composite of $A:f\tfl g$ and $B:h\tfl k$, with $\rep{g}=\rep{h}$, is defined by
\[
A\cdot B \:=\: (\sigma_g^- \star_1 A) \star_2 (\sigma_h^- \star_1 B).
\]
\item The identity of a $f:n\dfl 1$ is $\rep{\id}_f$.
\end{itemize}
If $f:n\dfl 1$ is an algebraic $2$-cell of $\Pr$, then the functor $\Pr(f):\Pr(n)\dfl\Pr(1)$ is defined by
\[
\Pr(f) \left( x_1,\dots, x_n \right) \:=\: (x_1\star_0\dots\star_0 x_n)\star_1 f. 
\]
Note that, using the naturality relations for $2$-cells of $\P$, we have:
\[
\Pr(f) \left( x_1,\dots, x_n \right) \:\approx\: (x_{\sigma_f(1)}\star_0\dots\star_0 x_{\sigma_f(n)}) \star_1 \rep{f}.
\]
If $A:f\tfl g:n\dfl A$ is an algebraic $3$-cell of $\Pr$, then the component at $(x_1,\dots,x_n)$ of the natural transformation~$\Pr(A)$ is given by 
\[
\Pr(A)_{(x_1,\dots,x_n)} \:=\: (x_1\star_0\dots\star_0 x_n) \star_1 A
\]
with source  
\[
(x_1\star_0\dots\star_0 x_n)\star_1 f \:\approx\: (x_{\sigma_f(1)}\star_0\dots\star_0 x_{\sigma_f(n)}) \star_1 \rep{f}
\]
and target
\[
(x_1\star_0\dots\star_0 x_n)\star_1 g \:\approx\: (x_{\sigma_g(1)}\star_0\dots\star_0 x_{\sigma_g(n)}) \star_1 \rep{g}.
\]

\begin{theorem}
\label{theoremGeneCoherenceProblemAlg2Prop}
Let $\P$ be an algebraic track $3$-prop and let $(A,B)$ be a $3$-sphere of $\P$. Then $A=B$ if and only if $\Pr(A)=\Pr(B)$. 
\end{theorem}

\begin{proof}
Let us assume that $A,B:f\tfl g:m\dfl n$ are such that $\Pr(A)=\Pr(B)$. Then we have, by definition of $\Pr$, for every algebraic $2$-cells $x_1$, $\dots$, $x_m$ of $\P$: 
\[
(x_1\star_0\dots\star_0 x_m) \star_1 A \:\approx\: (x_1\star_0\dots\star_0 x_m) \star_1 B.
\]
In particular, we take $\id_1$ for each $x_i$ to get $A\approx B$. Since $A$ and $B$ have the same source and the same target, we must have $A=B$.
\end{proof}

\subsection{The coherence theorem for braided monoidal categories}

\deftwocell[braid1]{b1 : 2 -> 2}
\deftwocell[braid2]{b2 : 2 -> 2}

\newcommand{\Braid}{\mathbf{Brd}}

\subsubsection{The \pdf{2}-pro of braids}

We define the $2$-pro of braids as the $2$-pro denoted by $\Braid$ and presented by the $2$-polygraph with two $2$-cells \twocell{b1} and \twocell{b2} and the following three $3$-cells
\[
\twocell{b1 *1 b2} \:\tfl\: \twocell{2}
\qquad\qquad
\twocell{b2 *1 b1} \:\tfl\: \twocell{2}
\qquad\qquad
\twocell{(b1 *0 1) *1 (1 *0 b1) *1 (b1 *0 1)} \:\tfl\: \twocell{(1 *0 b1) *1 (b1 *0 1) *1 (1 *0 b1)}
\]
In particular, those $3$-cells also generate the five following equalities in $\Braid$:
\[
\twocell{(b1 *0 1) *1 (1 *0 b1) *1 (b2 *0 1)}
	\:=\: \twocell{(1 *0 b2) *1 (b1 *0 1) *1 (1 *0 b1)}
\qquad\qquad
\twocell{(b1 *0 1) *1 (1 *0 b2) *1 (b2 *0 1)}
	\:=\: \twocell{(1 *0 b2) *1 (b2 *0 1) *1 (1 *0 b1)}
\qquad\qquad
\twocell{(b2 *0 1) *1 (1 *0 b1) *1 (b1 *0 1)}
	\:=\: \twocell{(1 *0 b1) *1 (b1 *0 1) *1 (1 *0 b2)}
\]
and
\[
\twocell{(b2 *0 1) *1 (1 *0 b2) *1 (b1 *0 1)}
	\:=\: \twocell{(1 *0 b1) *1 (b2 *0 1) *1 (1 *0 b2)}
\qquad\qquad
\twocell{(b2 *0 1) *1 (1 *0 b2) *1 (b2 *0 1)}
	\:=\: \twocell{(1 *0 b2) *1 (b2 *0 1) *1 (1 *0 b2)}
\]
The opposite $2$-pro $\Braid^o$ is the $2$-pro $\Braid$ with composition $\star_1$ reversed.

\begin{proposition}
The underlying category $\Br(1)$ of the initial algebra $\Br$ of $\Brd\Cat/\Mon\Cat$ is isomorphic to the $2$-pro $\Braid^o$.  
\end{proposition}

\begin{proof}
We note that, in the quotient track $3$-prop $\Brd\Cat/\Mon\Cat$, there is exactly one purely algebraic $2$-cell for each natural number~$n$. In particular, for $n=0$ and $n=1$, those are $\id_0$ and $\id_1$, respectively, for $n=2$, that is \twocell{mu} and, for $n\geq 3$, that is the equivalence class of any algebraic $2$-cell of $\Brd\Cat$ that contains exactly $(n-1)$ copies of \twocell{mu}.

Thus, the underlying category $\Br(1)$ of the initial $\Brd\Cat/\Mon\Cat$-algebra $\Br$ has the natural numbers as objects. Moreover, it is equipped with a structure of $2$-pro by the product $\tens$ defined by
\[
m\tens n \:=\: m+n 
\qquad\qquad\text{and}\qquad\qquad
A \tens B \:=\: (A\star_0 B) \star_1 \twocell{mu}
\]
Graphically, if $A=\twocell{dots *1 phi}$ and $B=\twocell{dots *1 psi}$, this product is written:
\[
\twocell{dots *1 phi} \tens\twocell{dots *1 psi} \:=\: \twocell{(dots *0 dots) *1 (phi *0 psi) *1 mu}.
\]
Let us define a morphism $\Phi:\Braid^o\fl\Br(1)$ of $2$-pros. On generating $2$-cells, we define
\[
\Phi(\twocell{b1}) \:=\: \twocell{beta}
\qquad\qquad\text{and}\qquad\qquad
\Phi(\twocell{b2}) \:=\: \twocell{beta}^-
\]
Let us prove that this induces a morphism of $2$-pros by checking that this is compatible with the generating $3$-cells of $\Braid$. For the first $3$-cell, we have:
\begin{align*}
\Phi\left(\twocell{b1 *1 b2}\right) 
	\:&=\: \Phi(\twocell{b2})\cdot\Phi(\twocell{b1}) 
\\
	\:&\approx\: \twocell{beta}^- \star_2 \twocell{beta} 
\\
	\:&=\: \twocell{mu} 
\\
	\:&=\: \Phi(\id_2).
\end{align*}
We prove, in a similar way, the relation $\Phi\left(\twocell{b2 *1 b1}\right) \approx \Phi\left(\id_2\right)$. For the last $3$-cell, we compute, on the one hand: 
\begin{align*}
\Phi\left(\twocell{(b1 *0 1) *1 (1 *0 b1) *1 (b1 *0 1)}\right)
	\:&=\: (\Phi(\twocell{b1}) \tens 1) 
		\cdot (1 \tens\Phi(\twocell{b1}))
		\cdot (\Phi(\twocell{b1}) \tens 1)
\\
	\:&\approx\: \twocell{(tau *0 1) *1 (1 *0 tau) *1 (beta *0 1) *1 mu}
		\star_2 \twocell{(tau *0 1) *1 (1 *0 beta) *1 mu}
		\star_2 \twocell{(beta *0 1) *1 mu}
\\
	\:&=\: \twocell{(tau *0 1) *1 (mu *0 1) *1 beta}
        \star_2 \twocell{(beta *0 1) *1 mu} 
\\
	\:&=\: \twocell{(beta *0 1) *1 beta}
\end{align*}
We have used the relation induced by the $4$-cell $\twocell{daleth}_1$ for the third equality and the exchange relation between $\star_1$ and $\star_2$ for the last equality. On the other hand, using the same properties, we get:
\begin{align*}
\Phi\left(\twocell{(1 *0 b1) *1 (b1 *0 1) *1 (1 *0 b1)}\right)
	\:&=\: (1\tens\Phi(\twocell{b1})) 
		\cdot (\Phi(\twocell{b1})\tens 1)
		\cdot (1\tens\Phi(\twocell{b1}))
\\
	\:&\approx\: \twocell{(1 *0 tau) *1 (tau *0 1) *1 (1 *0 beta) *1 mu}
		\star_2 \twocell{(1 *0 tau) *1 (beta *0 1) *1 mu}
		\star_2 \twocell{(1 *0 beta) *1 mu} 
\\
	\:&=\: \twocell{(beta *0 1) *1 tau *1 mu} \star_2 \twocell{(mu *0 1) *1 beta}
\\
	\:&=\: \twocell{(beta *0 1) *1 beta}
\end{align*}
Conversely, let us define a morphism $\Psi:\Br(1)\fl\Braid^o$ of $2$-pros. Using the exchange relation between $\star_1$ and $\star_2$, one can write any algebraic $3$-cell $A$ of $\Brd\Cat/\Mon\Cat$ as a composite
\begin{equation}
\label{decomposition}
A \:=\: A_1 \star_2 \cdots\star_2 A_k,
\end{equation}
where each $A_i$ is an algebraic $3$-cell of $\Brd\Cat/\Mon\Cat$ that contains exactly one generating $3$-cell, \ie, exactly one copy of either $\twocell{beta}$ or $\twocell{beta}^-$. Moreover, this decomposition is unique up to the inverse relations and the exchange relations between $\star_0$ and $\star_2$ and between $\star_1$ and $\star_2$.

The $4$-cells $\twocell{daleth}_1$ and $\twocell{daleth}_2$ generate the following relations in $\Brd\Cat/\Mon\Cat$
\begin{equation}
\label{relations1}
\twocell{(mu *0 1) *1 beta} 
	\:=\: \twocell{(1 *0 tau) *1 (beta *0 1) *1 mu} 
		\star_2 \twocell{(1 *0 beta) *1 mu} 
\qquad\text{and}\qquad
\left(\twocell{(1 *0 mu) *1 beta}\right) ^-
	\:=\: \left(\twocell{(beta *0 1) *1 mu}\right) ^-
		\star_2 \left(\twocell{(tau *0 1) *1 (1 *0 beta) *1 mu}\right) ^-
\end{equation}
which, in turn, using the inverse relations, induce
\begin{equation}
\label{relations2}
\left(\twocell{(mu *0 1) *1 beta}\right)^-
	\:=\: \left(\twocell{(1 *0 beta) *1 mu}\right)^-
		\star_2 \left(\twocell{(1 *0 tau) *1 (beta *0 1) *1 mu}\right)^-
\qquad\text{and}\qquad
\twocell{(1 *0 mu) *1 beta}
	\:=\: \twocell{(tau *0 1) *1 (1 *0 beta) *1 mu}
		\star_2 \twocell{(beta *0 1) *1 mu} 
\end{equation}
Those four relations have several consequences. The first one is that, in the decomposition~\eqref{decomposition}, we can assume that each $A_i$ has shape  
\deftwocell[black]{mu5:5 -> 1}
\[
m\tens\twocell{beta}^{\epsilon}\tens n \:=\: \left(\twocell{(dots *0 beta *0 dots) *1 mu5}\right)^{\epsilon}
\]
with $\epsilon$ in $\ens{-,+}$. In other terms, the $2$-pro $\Br$ admits $\twocell{beta}$ and $\twocell{beta}^-$ as generators. We define a morphism $\Psi:\Br(1)\fl\Braid^o$ of $2$-pros by 
\[
\Psi(\twocell{beta}) \:=\: \twocell{b1}
\qquad\qquad\text{and}\qquad\qquad
\Psi(\twocell{beta}^-) \:=\: \twocell{b2}
\]
This morphism is well-defined if and only if it is compatible with the inverse relations and the exchange relations of $\Brd\Cat/\Mon\Cat$. For the inverse relations, we use the fact that $\twocell{b2}$ is the inverse of $\twocell{b1}$ in the $2$-pro of braids. 

For the exchange relations between $\star_0$ and $\star_2$, we use the relations~\eqref{relations1} and~\eqref{relations2} to deduce that they are generated by the four relations
\deftwocell[black]{mu4:4 -> 1}
\[
\left(\twocell{(beta *0 dots *0 (tau *1 mu)) *1 mu4}\right)^{\epsilon_1}
	\star_2 \left(\twocell{(mu *0 dots *0 beta) *1 mu4}\right)^{\epsilon_2}
\:=\:
\left(\twocell{((tau *1 mu) *0 dots *0 beta) *1 mu4}\right)^{\epsilon_2}
	\star_2 \left(\twocell{(beta *0 dots *0 mu) *1 mu4}\right)^{\epsilon_1}
\]
where $\epsilon_1$ and $\epsilon_2$ range over $\ens{-,+}$. We check that, for each one, $\Psi$ sends both sides to the same braid. For example, in the case $\epsilon_1=\epsilon_2=+$, we get
\[
\Psi\left(\twocell{(beta *0 dots *0 (tau *1 mu)) *1 mu4}
	\star_2 \twocell{(mu *0 dots *0 beta) *1 mu4} \right)
	\:=\: \twocell{(4 *0 b1) *1 (2 *0 dots *0 2) *1 (b1 *0 4)} 
\]
and
\[
\Psi\left(\twocell{((tau *1 mu) *0 dots *0 beta) *1 mu4}
		\star_2 \twocell{(beta *0 dots *0 mu) *1 mu4}\right)
	\:=\: \twocell{(b1 *0 4) *1 (2 *0 dots *0 2) *1 (4 *0 b1)} 
\]
The relations~\eqref{relations1} and~\eqref{relations2} also induce that the exchange relations between $\star_1$ and $\star_2$ in $\Brd\Cat/\Mon\Cat$ are generated by the eight relations
\[
\left(\twocell{(tau *0 1) *1 (mu *0 1) *1 beta}\right)^{\epsilon_1} 
	\star_2 \left(\twocell{(beta *0 1) *1 mu}\right)^{\epsilon_2}
\:=\: 
\left(\twocell{(beta *0 1) *1 tau *1 mu}\right)^{\epsilon_2} 
	\star_2 \left(\twocell{(mu *0 1) *1 beta}\right)^{\epsilon_1}
\]
and
\[
\left(\twocell{(1 *0 tau) *1 (1 *0 mu) *1 beta}\right)^{\epsilon_1} 
	\star_2 \left(\twocell{(1 *0 beta) *1 mu}\right)^{\epsilon_2}
\:=\: 
\left(\twocell{(1 *0 beta) *1 tau *1 mu}\right)^{\epsilon_2} 
	\star_2 \left(\twocell{(1 *0 mu) *1 beta}\right)^{\epsilon_1}
\]
where $\epsilon_1$ and $\epsilon_2$ range over $\ens{-,+}$. We check that $\Psi$ is compatible with them. For example, in the case of the first relation, with $\epsilon_1=\epsilon_2=+$, we get:
\[
\Psi\left(\twocell{(tau *0 1) *1 (mu *0 1) *1 beta}
	\star_2 \twocell{(beta *0 1) *1 mu}\right)
	\:=\: \twocell{(b1 *0 1) *1 (1 *0 b1) *1 (b1 *0 1)}
\qquad\qquad\text{and}\qquad\qquad
\Psi\left(\twocell{(beta *0 1) *1 tau *1 mu}
	\star_2 \twocell{(mu *0 1) *1 beta}\right)
	\:=\: \twocell{(1 *0 b1) *1 (b1 *0 1) *1 (1 *0 b1)} 
\qedhere
\]
\end{proof}

\noindent For a $3$-cell $A$ of $\Brd\Cat/\Mon\Cat$, we identify the natural transformation $\Br(A)$ to its component at $(1,\dots,1)$, hence to $A$ itself and, using the isomorphism $\Psi:\Br(1)\fl\Braid^o$, to a braid on $n$ strands. By extension, if $A$ is a $3$-cell of $\Brd\Cat$, we denote by $\Br(A)$ the braid associated to its image in the quotient $\Brd\Cat/\Mon\Cat$.

\begin{theorem}[Coherence theorem for braided monoidal categories,~\cite{JoyalStreet93}]
Let $(A,B)$ be a $3$-sphere of $\Brd\Cat$. Then $A=B$ if and only if the braids $\Br(A)$ and $\Br(B)$ are equal. 
\end{theorem}

\begin{small}

\renewcommand{\refname}{\textbf{References}}
\providecommand{\bysame}{\leavevmode\hbox to3em{\hrulefill}\thinspace}

\vfill
\noindent
\hbox to7em{\hrulefill}

\bigskip
\noindent
\textsc{Yves Guiraud} \\
\textsc{INRIA, Institut Camille Jordan}, CNRS, Université de Lyon, Université Lyon~1 \\
Bâtiment Braconnier, 43 boulevard du 11 novembre 1918, 69622 Villeurbanne Cedex, France \\
E-mail address: \url{guiraud@math.univ-lyon1.fr}

\bigskip
\noindent
\textsc{Philippe Malbos} \\
\textsc{Institut Camille Jordan}, CNRS, Université de Lyon, Université Lyon~1 \\
Bâtiment Braconnier, 43 boulevard du 11 novembre 1918, 69622 Villeurbanne Cedex, France \\
E-mail address: \url{malbos@math.univ-lyon1.fr}

\end{small}

\end{document}